\def\ARTICLECLS{} 
  \ifx\ARTICLECLS\undefined
    \ifx\FULLSIZE\undefined
\documentclass{amsproc} 
    \else
\documentclass[12pt]{amsproc} 
    \fi
  \else
\documentclass[12pt,leqno]{article} 
  \fi
  \ifx\ARTICLECLS\undefined
\usepackage{aspmproc}
  \else
\makeatletter
\renewcommand\section{\@startsection {section}{1}{\z@}%
                                   {-3.5ex \@plus -1ex \@minus -.2ex}%
                                   {2.3ex \@plus.2ex}%
                                   {\normalfont\large\bfseries}}
\renewcommand\subsection{\@startsection{subsection}{2}{\z@}%
                                     {-3.25ex\@plus -1ex \@minus -.2ex}%
                                     {1.5ex \@plus .2ex}%
                                     {\normalfont\normalsize\bfseries}}
\makeatother
\usepackage{amsmath}
  \fi
\usepackage{amssymb}
  \ifx\SHOWKEYS\undefined\else
\usepackage{showkeys}
  \fi
  \ifx\ARTICLECLS\undefined
    \ifx\FULLSIZE\undefined\else 
    \fi
  \else
  \ifx\SHOWKEYS\undefined
  \else
  \fi
  \fi
\allowdisplaybreaks
\newcommand\A{{\mathcal A}}

\renewcommand\AA{{\mathbb A}}

\newcommand\ab{\allowbreak}
\newcommand\ad{\mathop{\mathrm{ad}}\nolimits}

\newcommand\av{\alpha^\vee}
\newcommand\B{{\mathcal B}}
\renewcommand\b{{\mathfrak b}}
\newcommand\bra{\langle}
\newcommand\C{{\mathbb C}} 
\newcommand\cl{{\mathrm{cl}}}
\newcommand\conj{\gamma}
\renewcommand\d{\partial}
\newcommand\End{\mathop{\mathrm{End}}\nolimits}
\newcommand\eps{\varepsilon}
\newcommand\ev{\mathop{\mathrm{ev}}\nolimits}
\newcommand\F{{\mathbb F}} 
\newcommand\g{{\mathfrak g}}
\newcommand\gr{\mathop{\mathrm{gr}}\nolimits}
\newcommand\h{{\mathfrak h}}
\newcommand\Hom{\mathop{\mathrm{Hom}}\nolimits}
\newcommand\id{\mathop{\mathrm{id}}\nolimits}

\newcommand\isom{\cong}
\newcommand\isomto{\overset\sim\rightarrow}

\newcommand\K{{\mathcal K}}

\newcommand\ket{\rangle}
\newcommand\lie{\mathrm}

\newcommand\N{{\mathcal N}}
\newcommand\n{{\mathfrak n}}
\newcommand\ord{\mathop{\mathrm{ord}}\nolimits}

\newcommand\PII{{\mathrm P}_{\mathrm{II}}}
\newcommand\PIII{{\mathrm P}_{\mathrm{III}}}
\newcommand\PIV{{\mathrm P}_{\mathrm{IV}}}
\newcommand\PV{{\mathrm P}_{\mathrm{V}}}
\newcommand\PVI{{\mathrm P}_{\mathrm{VI}}}
\newcommand\Q{{\mathbb Q}} 
\makeatletter\newcommand\qbinom{\genfrac[]\z@{}}\makeatother
\newcommand\rank{\mathop{\mathrm{rank}}\nolimits}
\renewcommand\setminus{\smallsetminus}

\newcommand\ts{{\tilde s}}

\newcommand\tw{{\tilde w}}
\newcommand\transpose[1]{\,{\vphantom{#1}}^t#1}
\newcommand\Spec{\mathop{\mathrm{Spec}}\nolimits}
\newcommand\Z{{\mathbb Z}} 
\newtheorem{theorem}{Theorem}

\newtheorem{lemma}[theorem]{Lemma}
\newtheorem{cor}[theorem]{Corollary}

\newtheorem{example}[theorem]{Example}
\newtheorem{definition}[theorem]{Definition}
\newtheorem{remark}[theorem]{Remark}

\numberwithin{theorem}{section}
\numberwithin{equation}{section}
\numberwithin{figure}{section}
\numberwithin{table}{section}

\newcommand\enumiref[1]{{\upshape(\ref{#1})}}
\newcommand\secref[1]{Section \ref{#1}}
\newcommand\theoremref[1]{Theorem \ref{#1}}

\newcommand\lemmaref[1]{Lemma \ref{#1}}
\newcommand\corref[1]{Corollary \ref{#1}}

\newcommand\exampleref[1]{Example \ref{#1}}

\newcommand\remarkref[1]{Remark \ref{#1}}

\newcommand{\BOX}{%
  \ifmmode\else\leavevmode\unskip\penalty9999\hbox{}\nobreak\hfill\fi
  \quad\hbox{\BOXSYMBOL}}
\makeatletter
\def\BOXSYMBOL{\RIfM@\bgroup\else$\bgroup\aftergroup$\fi
  \vcenter{\hrule\hbox{\vrule height.85em\kern.6em\vrule}\hrule}\egroup}
\makeatother
  \ifx\undefined\ARTICLECLS
  \else
\makeatletter
\newenvironment{proof}[1][\proofname]{\par
  \normalfont
  \topsep6\p@\@plus6\p@ \trivlist
  \item[\hskip\labelsep{\bfseries #1}\@addpunct{\bfseries.}]\ignorespaces
}{%
  \qed\endtrivlist
}
\newcommand\proofname{Proof}
\newcommand\qed{\BOX}
\makeatother
  \fi
  \ifx\ARTICLECLS\undefined
\title
[Quantization of Weyl group symmetries of Painlev\'e systems]
{Quantum groups and quantization of \\ Weyl group symmetries of Painlev\'e systems}
\author[G.~Kuroki]{Gen Kuroki}
\address{Mathematical Institute\\Tohoku University\\980-0814 Sendai\\Japan}
\rcvdate{September 30, 2007}
\rvsdate{September 30, 2007}
\dedicatory{Dedicated to Professor Akihiro Tsuchiya on his retirement}
  \else
\title
{\bf Quantum groups and quantization of \\ Weyl group symmetries of Painlev\'e systems}
\author{Gen Kuroki}
\date{August 19, 2008}
  \fi
\begin{document}

  \ifx\ARTICLECLS\undefined\else
\maketitle
  \fi

\begin{abstract}
 We shall construct
 the quantized $q$-analogues of the birational Weyl group actions
 arising from nilpotent Poisson algebras, 
 which are conceptual generalizations, proposed by Noumi and Yamada,
 of the B\"acklund transformations for Painlev\'e equations.
 Consider a quotient Ore domain of the lower nilpotent part
 of a quantized universal enveloping algebra 
 for any symmetrizable generalized Cartan matrix.
 Then non-integral powers of the image of the Chevalley generators 
 generate the quantized $q$-analogue of the birational Weyl group action.
 Using the same method, we shall reconstruct 
 the quantized B\"acklund transformations 
 of $q$-Painlev\'e equations constructed by Hasegawa.
 We shall also prove that
 any subquotient integral domain 
 of a quantized universal enveloping algebra
 of finite or affine type is an Ore domain.
\end{abstract}

  \ifx\ARTICLECLS\undefined
\maketitle
  \fi

\tableofcontents
\setcounter{section}{-1}

\section{Introduction}
\label{sec:intro}

The main theme of this article is 
a representation theoretic method for 
quantizing discrete symmetries 
of Painlev\'e equations and isomonodromic deformations.
Such discrete symmetries, often called B\"acklund transformations, 
play central roles in the theory of Painlev\'e equations and
isomonodromic deformations.
We could expect that this is the case in quantum settings.

In this article the term ``quantization'' means {\em canonical quantization}
that replaces commutative Poisson algebras $\A^\cl$ 
with non-commutative algebras $\A$ and 
Poisson algebra automorphisms of $\A^\cl$ with algebra automorphisms of $\A$.
If $\A^\cl$ is an integral domain and 
$Q(\A^\cl)$ denotes the field of fractions of $\A^\cl$, 
then a birational action of a group $G$ on $\Spec\A^\cl$ 
is identified with an algebra automorphism action of $G$ on $Q(\A^\cl)$.
Therefore 
if a classical symmetry is represented by a birational action of $G$ 
on a Poisson integral affine scheme, 
then its quantization should be an algebra automorphism action
of $G$ on a non-commutative skew field.

Note that $q$-difference analogue ($q$-analogue for short) 
or $q$-difference deformation ($q$-deformation) does not 
always mean quantization. 
In this article we shall deal with 
four types of classical and quantum systems, 
ordinary differential versions of classical systems and their quantizations, 
and $q$-analogues of classical systems and their quantizations.

Sections \ref{sec:QUEA} and \ref{sec:localizations}
are devoted to summarizing preparatory results on 
quantized universal enveloping algebras
and localizations of non-commutative rings.
We shall show that 
a quantized universal enveloping algebra
of finite or affine type is an Ore domain
(\theoremref{theorem:affine-tempered}). 
The $q=1$ cases were treated in 
\cite{rochacaridi-wallach-ore} and \cite{berman-cox-toroidal}.
Moreover we shall show that their subquotient integral domains 
are also Ore domains
(\corref{cor:subquot-affine-tempered}).
In \secref{sec:justifying-non-integral-powers},
we shall explain how to justify non-integral powers in 
fields of fractions
along the lines of the work \cite{iohara-malikov-GKconj}
by Iohara and Malikov.


\subsection{Quantized $q$-analogue of birational Weyl group actions}
\label{sec:intro-NY}

In \secref{sec:NY}, we shall construct 
a quantized $q$-analogue of the birational Weyl group action
arising from a nilpotent Poisson algebra.

A series of works by Okamoto 
\cite{okamoto-I,
okamoto-II, 
okamoto-III,
okamoto-IV}
showed that the Painlev\'e equations 
$\PII$, $\PIII$, $\PIV$, $\PV$, and $\PVI$ 
have affine Weyl group symmetries of type
$A^{(1)}_1$, $C^{(1)}_2$, $A^{(1)}_2$, $A^{(1)}_3$, and $D^{(1)}_4$
respectively.
Each of the affine Weyl groups birationally acts on 
dependent variables and parameters of the corresponding Painlev\'e equation.
Its birational actions preserve
the continuous flow generated by the Painlev\'e equation
and are called B\"acklund transformations.

In \cite{noumi-yamada-birataction1}, Noumi and Yamada generalized 
the birational affine Weyl group actions 
of type $A^{(1)}_2$, $A^{(1)}_3$, and $D^{(1)}_4$ 
(the cases of $\PIV$, $\PV$, and $\PVI$)
to the Weyl group associated to an arbitrary generalized Cartan matrix
(GCM for short).
If the size of GCM is equal to $m$, then the Weyl group birationally acts on 
a $2m$-dimensional space of $m$ dependent variables and $m$ parameters.
For each $m\geqq 2$, they also constructed a system of differential equations 
with $m$ dependent variables and $m$ parameters
on which the affine Weyl group of type $A^{(1)}_{m-1}$ acts 
as B\"acklund transformations \cite{noumi-yamada-higher}.
It is called a higher order Painlev\'e equation of type $A^{(1)}_{m-1}$.
Following this work, Nagoya constructed its quantization in \cite{nagoya}.

In \cite{noumi-yamada-birataction2},
Noumi and Yamada also proposed the further generalization 
of the birational Weyl group action.

Fix an arbitrary GCM.
Let $\g$ be the Kac-Moody algebra associated to the GCM,
$\h$ its Cartan subalgebra, and 
$\n_\pm$ its upper and lower parts
(\cite{kac}).
The Kostant-Kirillov Poisson bracket $\{\ ,\ \}$ makes 
the symmetric algebra $S(\n_-)$ Poisson 
and hence $\n_-^*=\Spec S(\n_-)$ is regarded as a Poisson scheme.
Let $J^\cl$ be an arbitrary Poisson prime ideal of $S(\n_-)$ and
denote by $\A_0^\cl=S(\n_-)/J^\cl$ the residue class ring modulo $J^\cl$.
Then $\A_0^\cl$ is a Poisson integral domain
and hence $\Spec\A_0^\cl$ is a Poisson integral subscheme of $\n_-^*$. 
Denote by $\A^\cl=\A_0^\cl\otimes S(\h)$ 
the tensor product algebra of $\A_0^\cl$ and $S(\h)$.
The Poisson structure of $\A_0^\cl$ uniquely extends to 
that of $\A^\cl$ so that $S(\h)$ is Poisson-central in $\A^\cl$.

In \cite{noumi-yamada-birataction2},
Noumi and Yamada constructed a birational Weyl group action
on $\Spec\A^\cl=\Spec\A_0^\cl\times\h$.
They called $\A_0^\cl$ a nilpotent Poisson algebra.
$\Spec \A_0^\cl$ and $\h$ are identified with 
the space of dependent variables 
and that of parameters respectively.

Let us explain the quantized $q$-analogue of the above setting.

Let $A=[a_{ij}]_{i,j\in I}$ be a symmetrizable GCM
symmetrized by a family $\{d_i\}_{i\in I}$ of positive rational numbers.
Denote by $\{\alpha_i\}_{i\in I}$ the set of simple roots 
and by $\{\av_i\}_{i\in I}$ the set of simple coroots.
Let $d$ be the least common denominator of $\{d_i\}_{i\in I}$.
Set the base field $\F$ by $\F=\Q(q^{1/d})$
and $q_i\in\F$ by $q_i=q^{d_i}$.  

Let $U_q$ be the quantized universal enveloping algebra 
of type $A$ over $\F$,
$U_q^0$ its Cartan subalgebra, and
$U_q^\pm$ its upper and lower parts
(\cite{lusztig}).
Let $J_q$ be an arbitrary completely prime ideal of $U_q^-$ and
denote by $\A_{q,0}=U_q^-/J_q$ the residue class ring modulo $J_q$.
Assume that $\A_{q,0}$ is an Ore domain.
For example, if $A$ is of finite or affine type,
then $\A_{q,0}$ is always an Ore domain.
See \corref{cor:subquot-affine-tempered}.
For the construction of examples for an arbitrary case,
see \secref{sec:truncated-q-Serre}.

Denote by $\A_q=\A_{q,0}\otimes U_q^0$ 
the tensor product algebra of $\A_{q,0}$ and $U_q^0$.
Then $\A_q$ is also an Ore domain.
Denote by $Q(\A_q)$ the skew field of fractions of $\A_q$.
Let $\{f_i\}_{i\in I}$ be the images in $\A_q$ 
of the lower Chevalley generators $\{F_i\}_{i\in I}$ of $U_q^-$.
In particular, $f_i$ ($i\in I$) satisfy the $q$-Serre relations.
Assume that $f_i\ne 0$ for all $i\in I$.

Denote by $\ts_i$ the action on $U_q^0$ of 
the simple reflection $s_i\in W$ for $i\in I$.
The action of $\ts_i$ naturally extends to the action on $Q(\A_q)$ 
so that $\ts_i$ trivially acts on $\A_{q,0}$.

In \secref{sec:NY}, we shall obtain the following results:
\begin{enumerate}
 \item For each $i\in I$, 
  the conjugation action $\conj(f_i^{\av_i})$ 
  of the non-integral power $f_i^{\av_i}$ on $Q(\A_q)$ formally given by 
  $\conj(f_i^{\av_i})x = f_i^{\av_i} x f_i^{-\av_i}$ for $x\in Q(\A_q)$ 
  is well-defined. (See \secref{sec:f_i^a}.)
 \item For each $i\in I$, define the operator $S_i$ acting on $Q(\A_q)$ by
  $S_i = \ts_i\circ\conj(f_i^{-\av_i}) = \conj(f_i^{\av_i})\circ\ts_i$.
  Then $S_i$ ($i\in I$) satisfy the defining relations of the Weyl group.
  In particular the braid relations of $S_i$ ($i\in I$) are derived from 
  the Verma relations of the Chevalley generators.
  (See the proof of \theoremref{theorem:quantized-NY}.)
\end{enumerate}
Thus we can construct a Weyl group action on $Q(\A_q)$,
which is the quantized $q$-analogue of the birational Weyl group action
arising from a nilpotent Poisson algebra.
For details, see \secref{sec:NY}.


\subsection{Quantized birational Weyl group actions of Hasegawa}
\label{sec:intro-KNY-Hasegawa}

In \secref{sec:KNY-Hasegawa},
we shall reconstruct 
the quantized ($q$-analogue of) birational Weyl group actions of Hasegawa
\cite{hasegawa}.

In \cite{kajiwara-noumi-yamada-qPIV}, 
Kajiwara, Noumi, and Yamada introduced a $q$-analogue of the fourth 
Painlev\'e equation $\PIV$.
It is called a $q$-Painlev\'e IV equation $q\PIV$.
They also constructed a birational action of 
the affine Weyl group of type $A^{(1)}_2$ 
preserving the discrete flow generated by $q\PIV$. 

Based on this work, in \cite{hasegawa}, 
Hasegawa constructed the quantized $q$-analogue of
the birational Weyl group action 
associated to an arbitrary symmetrizable GCM.

Apparently Hasegawa's quantized $q$-analogues are different from 
those explained in the preceding subsection.
But we can reconstruct the former by the same method as the latter.

We follow the notation and the assumptions 
on a GCM, roots, and coroots, etc.\ in the preceding subsection.
Assume 
that if $i\ne j$ and $a_{ij}\ne 0$, 
then $\epsilon_{ij}=\pm1$ and $\epsilon_{ji}=-\epsilon_{ij}$, 
and otherwise $\epsilon_{ij}=0$.

Let $\B_q$ be the associative algebra generated 
by $\{k_i^{\pm1},\,f_i\}_{i\in I}$ 
with following defining relations:
\begin{align*}
 &
 k_ik_i^{-1}=k_i^{-1}k_i=1, \quad
 k_ik_j=k_jk_i, 
 \\ &
 k_if_jk_i^{-1}=q_i^{-a_{ij}}f_j, \quad
 f_if_j=q_i^{-\epsilon_{ij}a_{ij}}f_jf_i.
\end{align*}
Note that the last formula 
is a sufficient condition for the $q$-Serre relation.
Define the algebra $U_q^0$ by
$U_q^0 = \F\left[\{a_i^{\pm1}\}_{i\in I}\right]$, 
the Laurent polynomial ring over $\F$ 
generated by $\{a_i^{\pm1}\}_{i\in I}$.
We identify $a_i$ with $q_i^{\av_i}$.
Define the algebra $\widetilde{\A}_q$ 
by $\widetilde{\A}_q=\B_q\otimes\B_q\otimes U_q^0$, 
the tensor product algebra of $\B_q$, $\B_q$, and $U_q^0$.
Then $\widetilde{\A}_q$ is an Ore domain
and hence can be embedded 
in the skew field of fractions $Q(\widetilde{\A}_q)$.

Define $f_{i1},f_{i2}\in\B_q\otimes\B_q$ by 
$f_{i1} = f_i\otimes 1$ and $f_{i2} = k_i^{-1}\otimes f_i$.
Note that $f_{i1}+f_{i2}$ is the image of the coproduct of 
a lower Chevalley generator in the quantized universal enveloping algebra.
Identify $f_{i1}\otimes 1, f_{i2}\otimes 1 , 1\otimes 1\otimes a_i \in \widetilde{\A}_q$ 
with $f_{i1},f_{i2}\in\B_q\otimes\B_q$, and $a_i\in U_q^0$ respectively.
Define $F_i\in Q(\widetilde{\A}_q)$ 
by $F_i = a_i^{-1}f_{i1}^{-1}f_{i2}$
for $i\in I$.
Do not confuse these $F_i$ with the lower Chevalley generators.
Let $\A_q$ be the subalgebra of $Q(\widetilde{\A}_q)$
generated by $\{F_i,a_i^{\pm1}\}_{i\in I}$.
Note that $a_i$ is central in the field of fractions $Q(\A_q)$ and 
\begin{equation*}
 F_i F_j = q_i^{-2\epsilon_{ij}a_{ij}} F_i F_i.
\end{equation*}
Therefore $\A_q$ is also an Ore domain.
We have $Q(\A_q)\subset Q(\widetilde{\A}_q)$.
Denote by $\ts_i$ the action on $U_q^0$ of 
the simple reflection $s_i\in W$ for $i\in I$.
Then we have $\ts_i(a_j)=a_ja_i^{-a_{ij}}$.
The action of $\ts_i$ extends to the action on $Q(\A_q)$ 
by $\ts_i(F_j)=a_i^{a_{ij}}F_j$.

In \secref{sec:KNY-Hasegawa}, we shall obtain the following results:
\begin{enumerate}
 \item For each $i\in I$, the conjugation action 
  $\conj\left((f_{i1}+f_{i2})^{\av_i}\right)$ 
  of the non-integral power $(f_{i1}+f_{i2})^{\av_i}$ 
  on $Q(\A_q)$ formally given by 
  \(\conj\left((f_{i1}+f_{i2})^{\av_i}\right)x 
   = (f_{i1}+f_{i2})^{\av_i} x (f_{i1}+f_{i2})^{-\av_i}\) for $x\in Q(\A_q)$ 
  is well-defined (\secref{sec:(f+f)^a}).
 \item For each $i\in I$, define the operator $S_i$ acting on $Q(\A_q)$ by \(
   S_i 
   = \ts_i\circ\conj\left((f_{i1}+f_{i2})^{-\av_i}\right)
   = \conj\left((f_{i1}+f_{i2})^{\av_i}\right)\circ\ts_i
  \).
  Then $S_i$ ($i\in I$) satisfy the defining relations of the Weyl group.
  In particular the braid relations of $S_i$ ($i\in I$) are derived from 
  the Verma relations of $\{f_{i1}+f_{i2}\}_{i\in I}$.
  (See \secref{sec:reconstruction-hasegawa}.)
\end{enumerate}
Thus we can construct a Weyl group action on $Q(\A_q)$,
which coincides with Hasegawa's quantized $q$-analogue of
the birational Weyl group action 
(\remarkref{remark:reconstruction-hasegawa}).
For details, see \secref{sec:KNY-Hasegawa}.


\section{Quantized universal enveloping algebras}
\label{sec:QUEA}

In this section, we shall summarize widely known results on 
quantized universal enveloping algebras.
We shall mainly follow Lusztig's book \cite{lusztig}.


\subsection{Symmetrizable GCM and root datum}
\label{sec:GCM-root}

A matrix $A=[a_{ij}]_{i\in I}$ with integer entries defined to be a 
{\em generalized Cartan matrix} ({\em GCM} for short) 
if it satisfies, for any $i,j\in I$,
(1) $a_{ii}=2$, 
(2) $a_{ij}\leqq 0$ if $i\ne j$, and 
(3) $a_{ji}=0$ if and only if $a_{ij}=0$.
Let $A$ be a GCM.
$A$ is called {\em indecomposable} if, for any $i\ne j$ in $I$,
there exists a sequence $i_0,i_1,\ldots,i_s\in I$
such that $i=i_0$, $a_{i_ki_{k+1}}\ne 0$ ($k=0,1,\ldots,s-1$), 
and $i_s=j$.
If there exists a family $\{d_i\}_{i\in I}$ of positive rational numbers 
such that $d_ia_{ij}=d_ja_{ji}$ for any $i,j\in I$, 
then $A$ is called {\em symmetrizable} 
and {\em symmetrized} by $\{d_i\}_{i\in I}$.
If $A$ is a GCM symmetrized by $\{d_i\}_{i\in I}$,
then the transpose $\transpose{A}$ is a GCM
symmetrized by $\{d_i^{-1}\}_{i\in I}$.

Let $A=[a_{ij}]_{i\in I}$ be a symmetrizable GCM 
symmetrized by $\{d_i\}_{i\in I}$.
We say that $A$ is {\em of finite type} 
(resp.\ {\em of affine type})
if its principal minors are positive
(resp.\ its proper principal minors are positive
and its determinant is equal to zero).
All GCM's of finite and affine type are classified explicitly.
For details, see Chapter 4 of \cite{kac}.

Let $X$ and $Y$ be finitely generated free $\Z$-modules
and $\bra\,,\ket:Y\times X\to \Z$ a perfect bilinear pairing.
($X$ can be identified with $\Hom_\Z(Y,\Z)$.)
Let $\{\av_i\}_{i\in I}$ and $\{\alpha_i\}_{i\in I}$ be
families of elements in $Y$ and $X$ respectively.
A {\em root datum} of type $A$ is defined to consist of
$(Y,\ab X,\ab\bra\,,\ket,\ab\{\av_i\}_{i\in I},\ab\{\alpha_i\}_{i\in I})$ 
satisfying $\bra\av_i,\alpha_j\ket = a_{ij}$ for any $i,j\in I$.
Then $\av_i$ and $\alpha_i$ are called a {\em coroot} and a {\em root} 
respectively.
The {\em dual root datum} of type $\transpose{A}$ is defined to be
$(X,\ab Y,\ab\bra\,,\ket,\ab\{\alpha_i\}_{i\in I},\ab\{\av_i\}_{i\in I})$.

If $\{\av_i\}_{i\in I}$ (resp.\ $\{\alpha_i\}_{i\in I}$) 
is linearly independent in $Y$ (resp.\ $X$), then
the root datum called {\em $Y$-regular} (resp.\ {\em $X$-regular})
and the submodule of $Y$ (resp.\ $X$) 
generated by $\{\av_i\}_{i\in I}$ (resp.\ $\{\alpha_i\}_{i\in I}$) 
is called a {\em coroot lattice} (resp.\ a {\em root lattice})
and denoted by $Q$ (resp.\ $Q^\vee$).
We set $X^+=\{\,\lambda\in X\mid
\bra\av_i,\lambda\ket\geqq 0\ \text{for all $i\in I$}\,\}$
and call its elements {\em dominant}.
We set $Q^+=\sum_{i\in I}\Z_{\geqq 0}\alpha_i$.


\subsection{Braid group and Weyl group}
\label{sec:braid-Weyl}

Let $A=[a_{ij}]_{i,j\in I}$ be a symmetrizable GCM.

The braid group $B(A)$ of type $A$ 
is the group generated by $\{s_i\}_{i\in I}$ with
the following defining relations: for any $i\ne j$ in $I$,
\begin{alignat*}{2}
 & s_is_j             = s_js_i & \quad & \text{if $a_{ij}a_{ji}=0$}, \\
 & s_is_js_i          = s_js_is_j & \quad & \text{if $a_{ij}a_{ji}=1$}, \\
 & s_is_js_is_j       = s_js_is_js_i & \quad & \text{if $a_{ij}a_{ji}=2$}, \\
 & s_is_js_is_js_is_j = s_js_is_js_is_js_i & \quad & \text{if $a_{ij}a_{ji}=3$}.
\end{alignat*}
These relations are called braid relations.

The Weyl group $W(A)$ of type $A$
is the group generated by $\{s_i\}_{i\in I}$ satisfying
the braid relations together with $s_i^2=1$ for all $i\in I$.
When $A$ is indecomposable, 
$W(A)$ is finite if and only if $A$ is of finite type.

Denote by $B$ the braid group of type $A$ 
and by $W$ the Weyl group of type $A$.

For $w\in W$, the length $\ell(w)$ of $w$ is the smallest integer $p\geqq 0$
such that there exists $i_1,\ldots,i_p\in I$ with $w=s_{i_1}\cdots s_{i_p}$.
Then $s_{i_1}\cdots s_{i_p}$ is called a reduced expression of $w$.

If $s_{i_1}\cdots s_{i_p}$ and $s_{i_1'}\cdots s_{i_p'}$ are 
reduced expressions of $w\in W$, 
then the equality $s_{i_1}\cdots s_{i_p} \ab=\ab s_{i_1'}\cdots s_{i_p'}$ 
holds in the braid group $B$. 
Therefore the mapping from $W$ to $B$ sending $w\in W$ to 
the element of $B$ represented by a reduced expression 
of $w$ is well-defined.

Let 
$(Y,\ab X,\ab \bra\,,\ab\ket,\ab\{\av_i\}_{i\in I},\ab\{\alpha_i\}_{i\in I})$
be a root datum of type $A$.
Then the Weyl group $W(A)$ acts on $Y$ and $X$ 
by $s_i(y)=y-\bra y,\alpha_i\ket\av_i$ for $y\in Y$
and $s_i(x)=x-\bra\av_i,x\ket\alpha_i$ for $x\in X$.
Moreover we have $\bra w(y),x\ket = \bra y,w^{-1}(x)\ket$ 
for $w\in W(A)$, $y\in Y$, and $x\in X$.

If $s_{i_1}\cdots s_{i_p}$ is a reduced expression in $W$,
then $s_{i_p}s_{i_{p-1}}\cdots s_{i_2}(\av_{i_1})\in \sum_{i\in I}\Z_{\geqq0}\av_i$
and $s_{i_p}s_{i_{p-1}}\cdots s_{i_2}(\alpha_{i_1})\in \sum_{i\in I}\Z_{\geqq0}\alpha_i$.

\begin{example}
\label{example:braid-relations}
\normalfont
 Assume that $i,j\in I$ and $i\ne j$.
 \begin{enumerate}
  \item If $(a_{ij},a_{ji})=(0,0)$, 
   then both sides of $s_is_j=s_js_i$ are reduced expressions and
   \begin{equation*}
    \begin{cases}
     1(\av_j)  = \av_j, \\
     s_j(\av_i)= \av_i, 
    \end{cases}
    \begin{cases}
     1(\av_i)  = \av_i, \\
     s_i(\av_j)= \av_j.
    \end{cases}
   \end{equation*}
  \item If $(a_{ij},a_{ji})=(-1,-1)$, 
   then both sides of $s_is_js_i=s_js_is_j$ are reduce expressions and
   \begin{equation*}
    \begin{cases}
     1(\av_i)      = \av_i, \\
     s_i(\av_j)    = \av_i+\av_j, \\
     s_is_j(\av_i) = \av_j, 
    \end{cases}
    \begin{cases}
     1(\av_j)      = \av_j, \\
     s_j(\av_i)    = \av_i+\av_j, \\
     s_js_i(\av_j) = \av_i.
    \end{cases}
   \end{equation*}
  \item If $(a_{ij},a_{ji})=(-1,-2)$, 
   then both sides of $s_is_js_is_j=s_js_is_js_i$ 
   are reduce expressions and
   \begin{equation*}
    \begin{cases}
     1(\av_j)         = \av_j, \\
     s_j(\av_i)       = \av_i+\av_j, \\
     s_js_i(\av_j)    = 2\av_i+\av_j, \\
     s_js_is_j(\av_i) = \av_i, 
    \end{cases}
    \begin{cases}
     1(\av_i)         = \av_i, \\
     s_i(\av_j)       = 2\av_i+\av_j, \\
     s_is_j(\av_i)    = \av_i+\av_j, \\
     s_is_js_i(\av_j) = \av_j.
    \end{cases}
   \end{equation*}
  \item If $(a_{ij},a_{ji})=(-1,-3)$, 
   then both sides of 
   $s_is_j\ab s_is_j\ab s_is_j\ab=\ab s_js_i\ab s_js_i\ab s_js_i$ 
   are reduce expressions and
   \begin{equation*}
    \begin{cases}
     1(\av_j)               = \av_j, \\
     s_j(\av_i)             = \av_i+\av_j, \\
     s_js_i(\av_j)          = 3\av_i+2\av_j, \\
     s_js_is_j(\av_i)       = 2\av_i+\av_j, \\
     s_js_is_js_i(\av_j)    = 3\av_i+\av_j, \\
     s_js_is_js_is_j(\av_i) = \av_i,
    \end{cases}
    \begin{cases}
     1(\av_i)               = \av_i, \\
     s_i(\av_j)             = 3\av_i+\av_j, \\
     s_is_j(\av_i)          = 2\av_i+\av_j, \\
     s_is_js_i(\av_j)       = 3\av_i+2\av_j, \\
     s_is_js_is_j(\av_j)    = \av_i+\av_j, \\
     s_is_js_is_js_i(\av_j) = \av_j.
    \end{cases}
   \end{equation*}
 \end{enumerate}
 These formulae shall be applied to the Verma relations.
 \BOX
\end{example}


\subsection{Kac-Moody algebra}
\label{sec:Kac-Moody}

For details of Kac-Moody algebras, see Kac's book \cite{kac}.

Let $A=[a_{ij}]_{i,j\in I}$ be 
a symmetrizable GCM and 
$(Y,\ab X,\ab\bra\,,\ket,\ab\{\av_i\}_{i\in I},\ab\{\alpha_i\}_{i\in I})$
a root datum of type $A$.
We set $\h=\C\otimes_\Z Y$ and 
identify $\h^*$ with $\C\otimes_\Z X$ by $\bra\,,\ket$.

The Kac-Moody (Lie) algebra $\g$ associated to the root datum is defined to be
the Lie algebra over $\C$
generated by $E_i$, $F_i$ ($i\in I$) and $H\in\h$
with following defining relations:
\begin{align*}
 &
 \text{$\h$ is an Abelian Lie subalgebra of $\g$};
 \\ &
 [H,E_i] =  \bra H,\alpha_i\ket E_i, \quad
 [H,F_i] = -\bra H,\alpha_i\ket F_i
 \quad\text{for $i\in I$, $H\in\h$};
 \\ &
 [E_i, F_j] = \delta_{ij}\av_i
 \quad\text{for $i,j\in I$};
 \\ &
 \ad(E_i)^{1-a_{ij}}E_j=0,
 \quad
 \ad(F_i)^{1-a_{ij}}F_j=0
 \quad\text{if $i\ne j$}.
\end{align*}
Here we set $\ad(X)Y=[X,Y]$, for example, $\ad(X)^3Y=[X,[X,[X,Y]]]$.
The last two relations are called {\em Serre relations}.

Denote by $\n_+$ (resp.\ $\n_-$)
the Lie subalgebra of $\g$
generated by $\{E_i\}_{i\in I}$ (resp.\ $\{F_i\}_{i\in I}$).
We call $\n_\pm$ the {\em upper and lower parts} of $\g$
and $E_i$ (resp.\ $F_i$) 
the Chevalley generators of $\n_+$ (resp. $\n_-$).
The Abelian subalgebra $\h$ is called the {\em Cartan subalgebra} of $\g$.
We have the triangular decomposition $\g=\n_-\oplus\h\oplus\n_+$ of $\g$.
Define the {\em upper and lower Borel subalgebras} $\b_\pm$ 
by $\b_\pm=\h\oplus\n_\pm$.

We can define the $\Z$-gradation of $\g$ by $\deg E_i=1$, $\deg F_i=-1$
($i\in I$) and $\deg H=0$ ($H\in\h$) 
and call it the {\em principal gradation} of $\g$.
Denote by $\g_k$ the degree-$k$ part of $\g$ for $k\in\Z$.
Then we have $\g = \bigoplus_{k\in\Z}\g_k$, $\dim\g_k<\infty$,
$\n_\pm=\bigoplus_{k>0}\g_{\pm k}$, 
and $\h=\g_0$.
The induced $\Z$-gradation $U(\g)=\bigoplus_{k\in\Z}U(\g)_k$ is
also called the principal gradation.
Define the principal gradations of $U(\n_\pm)$ 
by $U(\n_\pm)_{\pm k}=U(\g)_{\pm k}\cap U(\n_\pm)$ 
for $k\in\Z_{\geqq0}$.

Assume that the root datum is $Y$-regular and $X$-regular.
Denote by $\g'$ the derived Lie algebra $[\g,\g]$ of 
the Kac-Moody algebra $\g$.
Then $\g'$ is a finite dimensional Lie algebra
(resp.\ a central extension of a (possibly twisted) loop algebra 
of a finite dimensional simple Lie algebra)
if and only if
$A$ is of finite type (resp.\ of affine type).
Then $\g$ is called {\em a Kac-Moody algebra of finite type} 
(resp.\ {\em a Kac-Moody algebra of affine type} 
or {\em an affine Lie algebra} for short).
These lead to the following result.

\begin{lemma}
\label{lemma:dimU(n)k}
 Assume that the GCM $A$ is of finite or affine type.
 Then $\{\dim\g_k\}_{k\in\Z}$ is bounded,
 namely there exists a positive integer $N$ such that
 $\dim\g_k\leqq N$ for all $k\in\Z$.
 Define the positive integers $C^{(N)}_k$ ($k\in\Z_{\geqq0}$) by 
 $\left(\prod_{i=1}^\infty(1-t^i)\right)^{-N}
 =\sum_{k=0}^\infty C^{(N)}_k t^k$.
 Then we have $\dim U(\n_\pm)_{\pm k}\leqq C^{(N)}_k$ for $k\in\Z_{\geqq0}$.
 \BOX
\end{lemma}


\subsection{$q$-Binomial theorem}
\label{q:q-binom}

We define $q$-numbers, $q$-factorials, $q$-binomial coefficients, 
and $q$-shifted factorials as follows:
\begin{align*}
 &
 [x]_q = \frac{q^x-q^{-x}}{q-q^{-1}}
 \quad\text{for $n\in\Z$},
 \\ &
 [n]_q! = \prod_{k=1}^n [k]_q
 \quad\text{for $n\in\Z_{\geqq0}$},
 \\ &
 \qbinom{x}{k}_q
 = \frac{[x]_q[x-1]_q\cdots[x-k+1]_q}{[k]_q!} 
 \quad\text{for $k\in\Z_{\geqq0}$},
 \\ &
 (x)_{q,k} = (1+x)(1+q^2x)\cdots(1+q^{2(k-1)}x)
 \quad\text{for $k\in\Z_{\geqq0}$},
 \\ &
 (x)_{q,\infty} 
 = (1+x)(1+q^2x)(1+qx^4)\cdots 
 = \prod_{\mu=0}^\infty(1+q^{2\mu}x).
\end{align*}
Note that our $q$-shifted factorials are 
different from usual ones defined 
by $(x;q)_k=\prod_{\mu=0}^{k-1}(1-q^\mu x)$.
Then we can prove the following lemma by induction on $n$.

\begin{lemma}[$q$-binomial theorem]
\label{lemma:q-binom}
 Assume that $x,y$ are elements of an $\Q(q)$-algebra satisfying $yx=q^2xy$.
 Then, for $n=0,1,2,\ldots$,
 \begin{align*}
  (x+y)^n 
  &
  = \sum_{k=0}^n q^{k(n-k)} \qbinom{n}{k}_q x^k y^{n-k}
  \\ &
  = \sum_{k=0}^\infty q^{k(n-k)} \qbinom{n}{k}_q x^k y^{n-k}
  = \sum_{k=0}^\infty q^{k(n-k)} \qbinom{n}{k}_q x^{n-k} y^k.
 \end{align*}
 Moreover, if $x$ is invertible, then, for $n=0,1,2,\ldots$,
 \begin{equation*}
  (x+y)^n 
  = x^n (x^{-1}y)_{q,n} 
  = x^n \frac{(x^{-1}y)_{q,\infty}}{(q^{2n}x^{-1}y)_{q,\infty}}
  = \frac{(q^{-2n}x^{-1}y)_{q,\infty}}{(x^{-1}y)_{q,\infty}} x^n,
 \end{equation*}
 where the infinite products cancel out except finite factors.
 \BOX
\end{lemma}


\subsection{Quantized universal enveloping algebra}
\label{sec:def-QUEA}

Let $A=[a_{ij}]_{i,j\in I}$ be 
a symmetrizable GCM symmetrized by $\{d_i\}_{i\in I}$
and $(Y,\ab X,\ab\bra\,,\ket,\ab\{\av_i\}_{i\in I},\ab\{\alpha_i\}_{i\in I})$
a root datum of type $A$.
Let $d$ be the least common denominator of $\{d_i\}_{i\in I}$.
We set the base field $\F$ by $\F=\Q(q^{1/d})$
and $q_i\in\F$ by $q_i=q^{d_i}$.  
Then we have $d_i\av_i\in d^{-1}Y$
and we extend naturally
the perfect bilinear pairing $\bra\,,\ket:Y\times X\to\Z$ 
to $\bra\,,\ket:d^{-1}Y\times X\to d^{-1}\Z$.

Then the {\em quantized universal enveloping algebra} $U_q=U_q(\g)$
associated to the root datum
is defined to be the associative algebra over the base field $\F$
generated by $E_i$, $F_i$, $q^\lambda$ for $i\in I$ and $\lambda\in d^{-1}Y$
with following defining relations:
\begin{align*}
 &
 q^0 = 1, \quad 
 q^{\lambda+\mu}=q^\lambda q^\mu \quad
 \text{for $\lambda,\mu\in d^{-1}Y$};
 \\ &
 q^\lambda E_i q^{-\lambda} = q^{\bra\lambda,\alpha_i\ket} E_i, \quad
 q^\lambda F_i q^{-\lambda} = q^{-\bra\lambda,\alpha_i\ket} F_i 
 \quad\text{for $i\in I$, $\lambda\in d^{-1}Y$};
 \\ &
 E_iF_j - F_jE_i 
 = \delta_{ij}[\av_i]_{q_i}
 = \delta_{ij}\frac{K_i - K_i^{-1}}{q_i-q_i^{-1}}
 \quad\text{for $i,j\in I$};
 \\ &
 \sum_{k=0}^{1-a_{ij}} 
 (-1)^k \qbinom{1-a_{ij}}{k}_{q_i} E_i^k E_j E_i^{1-a_{ij}-k} = 0
 \quad\text{if $i\ne j$};
 \\ &
 \sum_{k=0}^{1-a_{ij}} 
 (-1)^k \qbinom{1-a_{ij}}{k}_{q_i} F_i^k F_j F_i^{1-a_{ij}-k} = 0
 \quad\text{if $i\ne j$},
\end{align*}
where we set $K_i = q_i^{\av_i} = q^{d_i\av_i}$.
The last two relations are called {\em $q$-Serre relations}.
(For this definition, see Corollary 33.1.5 in \cite{lusztig}.)
In particular, we have 
$K_iE_jK_i^{-1}=q_i^{a_{ij}}E_j$ and
$K_iF_jK_i^{-1}=q_i^{-a_{ij}}F_j$ for $i,j\in I$.
By induction on $m,n\in\Z_{\geqq0}$ we can prove the following formula:
\begin{align}
 &
 E_i^{(m)} F_i^{(n)} = 
 \sum_{k=0}^{\min\{m,n\}}
 F_i^{(n-k)} \qbinom{\av_i-(m-k)-(n-k)}{k}_{q_i} E_i^{(m-k)},
 \notag
 \\ &
 F_i^{(n)} E_i^{(m)} = 
 \sum_{k=0}^{\min\{m,n\}}
 E_i^{(m-k)} \qbinom{-\av_i-(m-k)-(n-k)}{k}_{q_i} F_i^{(n-k)},
 \label{eq:F^nE^m}
\end{align}
where we set $E_i^{(m)}=E_i^m/[m]_{q_i}!$ and $F_i^{(n)}=F_i^n/[n]_{q_i}!$.
(See Corollary 3.1.9 in \cite{lusztig}.)

Denote by $U_q^+=U_q(\n_+)$ (resp.\ $U_q^-=U_q(\n_-)$, $U_q^0=U_q(\h)$)
the subalgebra of $U_q=U_q(\g)$ 
generated by $\{E_i\}_{i\in I}$ 
(resp.\ $\{F_i\}_{i\in I}$, $\{q^\lambda\}_{\lambda\in d^{-1}Y}$).
We call $U_q^\pm$ the {\em upper and lower parts} of $U_q$
and $E_i$, $F_i$ the upper and lower {\em Chevalley generators} respectively.
The commutative subalgebra $U_q^0$ is called the {\em Cartan subalgebra} of $U_q$.
We have the triangular decomposition 
$U_q\isom U_q^-\otimes U_q^0\otimes U_q^+$ of $U_q$.
Note that $U_q^\pm$ are determined by the symmetrized GCM only
and their structure does not depend on the choice of a root datum.
Define the {\em upper and lower Borel subalgebras} $U_q(\b_\pm)$
to be the subalgebra of $U_q=U_q(\g)$ 
generated by $U_q^\pm=U_q(\n_\pm)$ and $U_q^0=U_q(\h)$.
Then we have 
$U_q(\b_\pm)\isom U_q^0\otimes U_q^\pm \isom U_q^\pm\otimes U_q^0$.
We say that $U_q$ is {\em of finite type} 
(resp.\ {\em of affine type} or {\em affine} for short)
if $A$ is of finite type (resp.\ of affine type).

We can define the $\Z$-gradation of $U_q$ by $\deg E_i=1$, $\deg F_i=-1$
($i\in I$) and $\deg q^\lambda=0$ ($\lambda\in d^{-1}Y$) 
and call it the {\em principal gradation} of $U_q$.
Denote by $(U_q)_k$ the degree-$k$ part of $U_q$ for $k\in\Z$.
Then we have $U_q = \bigoplus_{k\in\Z}U_k$.
Define the principal gradations of $U_q^\pm$ 
by $(U_q^\pm)_{\pm k}=(U_q)_{\pm k}\cap U_q^\pm$ 
for $k\in\Z_{\geqq0}$.
Then we have $U_q^\pm=\bigoplus_{k\geqq 0}(U_q^\pm)_{\pm k}$.

We can regard $U_q$ (resp.\ $U_q^\pm$, $U_q^0$) 
as a $q$-deformation of 
the universal enveloping algebra $U(\g)$, 
(resp.\ $U(\n_\pm)$, $U(\h)$)
of a Kac-Moody algebra $\g$
(resp. its upper and lower parts $\n_\pm$, its Cartan subalgebra $\h$).

Define the local ring $\AA_1$ by \(
 \AA_1
 =\{\,f(q^{1/d})\in\F=\Q(q^{1/d})\ab
  \mid f\ \text{is}\ab\ \text{regular}\ab\ \text{at}\ab\ q^{1/d}=1\,\}
 =\Q(q^{1/d})\ab\cap\ab\Q[[q^{1/d}-1]]
\).
We regard $\C$ as an algebra over $\AA_1$ by acting $q$ on $\C$ as $1$
and denote the algebra $\C$ over $\AA_1$ by $\C_1$.
Assume that $\{y_1,\ldots,y_M\}$ is a $\Z$-free basis of $d^{-1}Y$.
Set $(x)_q=(1-q^x)/(1-q)$.
Then we have the following results.
For the proof, see Sections 3.3 and 3.4 of \cite{hong-kang},
for example.

\begin{lemma}
\label{lemma:U_AA1}
 Let $U_{\AA_1}$ be the subalgebra of $U_q$ over $\AA_1$
 generated by $\{E_i,F_i\}_{i\in I}$
 and $\{(y_\mu)_q,q^{-y_\mu}\}_{\mu=1}^M$.
 Let $U_{\AA_1}^+$ (resp.\ $U_{\AA_1}^-$) be 
 the subalgebra of $U_q^+$ over $\AA_1$ (resp.\ $U_q^-$)
 generated by $\{E_i\}_{i\in I}$ (resp.\ $\{F_i\}_{i\in I}$)
 and $U_{\AA_1}^0$ 
 the subalgebra of $U_q^0$ over $\AA_1$
 generated by $\{q^{y_\mu},(y_\mu)_q\}_{\mu=1}^M$.
 Set $(U_{\AA_1}^\pm)_k = (U_q^\pm)_k\cap U_{\AA_1}$ for $k\in\Z$.
 \begin{enumerate}
 \item \label{enumi:U_AA1-1}
  The multiplication gives an isomorphism 
  $U_{\AA_1}^-\otimes_{\AA_1}U_{\AA_1}^0\otimes_{\AA_1}U_{\AA_1}^+\isomto U_{\AA_1}$
  of $\AA_1$-modules.
 \item \label{enumi:U_AA1-2}
  $(U_{\AA_1}^\pm)_k$ are free $\AA_1$-modules
  and $U_{\AA_1}^\pm=\bigoplus_{k=0}^\infty (U_{\AA_1}^\pm)_{\pm k}$.
 \item \label{enumi:U_AA1-3}
  \(
   U_{\AA_1}^0
   =\AA_1[(y_1)_q,\ldots,(y_M)_q,q^{-y_1},\ldots,q^{-y_M}]
  \) properly contains \(
   \AA_1[q^{\pm y_1},\ldots,q^{\pm y_M}]
  \).
 \item \label{enumi:U_AA1-4}
  $\F\otimes_{\AA_1} U_{\AA_1}=U_q$,
  $\F\otimes_{\AA_1}U_{\AA_1}^\pm=U_q^\pm$,
  $\F\otimes_{\AA_1}(U_{\AA_1}^\pm)_k=(U_q^\pm)_k$,
  and $\F\otimes_{\AA_1}U_{\AA_1}^0=U_q^0$.
 \item \label{enumi:U_AA1-5}
  $\C_1\otimes_{\AA_1} U_{\AA_1}=U(\g)$,
  $\C_1\otimes_{\AA_1}U_{\AA_1}^\pm=U(\n_\pm)$,
  $\C_1\otimes_{\AA_1}(U_{\AA_1}^\pm)_k=U(\n_\pm)_k$,
  and $\C_1\otimes_{\AA_1}U_{\AA_1}^0=U(\h)$.
 \item \label{enumi:U_AA1-6}
  \(
   \dim_\F (U_q^\pm)_{\pm k} 
   = \rank_{\AA_1} (U_{\AA_1}^\pm)_{\pm k}
   = \dim_\C U(\n_\pm)_{\pm k}
   \) for $k\in\Z_{\geqq0}$.
 \BOX
 \end{enumerate}
\end{lemma}

In particular, we obtain the following results.

\begin{lemma}
\label{lemma:dimUq(n)k}
 Let $A$ be a symmetrizable GCM.
 A quantized universal enveloping algebra $U_q$ of type $A$
 is always an integral domain.
 If $A$ is of finite or affine type, 
 then $\dim_\F(U_q^\pm)_k\leqq C^{(N)}_k$
 for $k\in\Z_{\geqq0}$,
 where $N$ and $C^{(N)}_k$ are given in \lemmaref{lemma:dimU(n)k}.
\end{lemma}

\begin{proof}
 Assume that $a,b\in U_q$ are non-zero.
 Then there exist $\lambda,\mu\in d^{-1}Y$ such that \(
 q^\lambda a, bq^\mu 
 \in U_q^-\otimes \F[q^{y_1},\ldots,q^{y_M}]\otimes U_q^+
 \).  
 Let $\{u^\pm_s\}_{s=0}^\infty$ be $\AA_1$-free bases of $U_{\AA_1}^\pm$.
 Set $u^0_\mu = (y_1)_q^{\mu_1}\cdots(y_M)_q^{\mu_M}$
 for $\mu=(\mu_1,\ldots,\mu_M)\in(\Z_{\geqq0})^n$.
 Then $\{u^0_\mu\}_{\mu\in(\Z_{\geqq0})^n}$ is
 an $\AA_1$-free basis of $\AA_1[(y_1)_q,\ldots,(y_M)_q]$.
 Note that $\F[q^{y_1},\ldots,q^{y_M}] = \F[(y_1)_q,\ldots,(y_M)_q]$. 
 Because of \lemmaref{lemma:U_AA1} 
 \enumiref{enumi:U_AA1-1} and \enumiref{enumi:U_AA1-4},
 we can uniquely write $q^\lambda a$ and $b q^\mu$
 in the following forms: \(
  q^\lambda a
  = \sum_{r,\mu,s} c_{r,\mu,s} u^-_r u^0_\mu u^+_s
 \), \(
  b q^\mu
  = \sum_{r,\mu,s} d_{r,\mu,s} u^-_r u^0_\mu u^+_s
 \) ($c_{r,\mu,s},d_{r,\mu,s}\in\F=\Q(q^{1/d})$),
 where only finitely many $c_{r,\mu,s}$ and $d_{r,\mu,s}$ are non-zero.
 Since any $c\in\F^\times$ is uniquely expressed as $c=(q-1)^{-k}\tilde{c}$
 with $k\in\Z$ and $\tilde{c}\in\AA_1^\times$ 
 (i.e.\ $\tilde{c}\in\AA_1$ and $\tilde{c}(1)\ne 0$),
 we set $\ord(c)=k$ and $\ord(0)=-\infty$.
 Setting $l=\max\{\ord(c_{s,\mu,s})\}$ and $m=\max\{\ord(d_{s,\mu,s})\}$,
 we have \(
  (q-1)^l q^\lambda a,
  (q-1)^m b q^\mu
  \in U_{\AA_1}^-\o_{\AA_1}\AA_1[(y_1)_q,\ldots,(y_M)_q]\o_{\AA_1}U_{\AA_1}^+
 \).
 Moreover their images in $\C_1\o_{\AA_1}U_{\AA_1}=U(\g)$ are non-zero
 and hence their product in $U(\g)$ is non-zero.
 Therefore $(q-1)^{m+l}q^\lambda ab q^\mu \ne 0$, namely $ab\ne 0$.
 This means that $U_q$ is an integral domain.
 The second statement immediately follows from 
 \lemmaref{lemma:U_AA1} \enumiref{enumi:U_AA1-6}.
\end{proof}

If the root datum is $Y$-regular and $X$-regular,
then the highest weight integral representations of $\g$ 
are deformed to those of $U_q$.
For details, see Chapter 33 of \cite{lusztig} 
and Section 3.4 of \cite{hong-kang}.

We can define 
the coproduct $\Delta:U_q\to U_q\otimes U_q$ (an algebra homomorphism),
the counit $\eps:U_q\to\F$ (an algebra homomorphism), and
the antipode $S:U_q\to U_q$ (an anti-algebra automorphism) by
\begin{align*}
 &
 \Delta(E_i) = E_i\otimes K_i + 1\otimes E_i
 \quad\text{for $i\in I$},
 \\ &
 \Delta(F_i) = F_i\otimes 1 + K_i^{-1}\otimes F_i
 \quad\text{for $i\in I$},
 \\ &
 \Delta(q^\lambda) = q^\lambda\otimes q^\lambda
 \quad\text{for $\lambda\in d^{-1}Y$},
 \\ &
 \eps(E_i)=0, \quad \eps(F_i)=0, \quad \eps(q^\lambda)=1
 \quad\text{for $i\in I$, $\lambda\in d^{-1}Y$},
 \\ &
 S(E_i)=-E_iK_i^{-1},\ S(F_i)=-K_iF_i,\ S(q^\lambda)=q^{-\lambda}
 \ \text{for $i\in I$, $\lambda\in d^{-1}Y$}.
\end{align*}
These give a Hopf algebra structure on $U_q$.

\begin{remark}
\label{remark:Hopf-algebra-structures}
\normalfont
 The above definition of a Hopf algebra structure on $U_q$ 
 is different from that in Lusztig's book \cite{lusztig}.
 Denote by $\Delta^L$, $\eps^L$, and $S^L$ 
 the coproduct, the counit, and the antipode of \cite{lusztig} respectively.
 We can uniquely define 
 the involutive algebra automorphism $\omega$ of $U_q$ by
 $\omega(E_i) = F_i$, $\omega(F_i) = E_i$, 
 $\omega(q^\lambda) = q^{-\lambda}$ for $i\in I$, $\lambda\in d^{-1}Y$.
 Then the Hopf algebra structure of \cite{lusztig} 
 is related to ours by 
 $\Delta^L=(\omega\otimes\omega)\circ\Delta\circ\omega$,
 $\eps^L=\eps\circ\omega$, and $S^L=\omega\circ S\circ\omega$.
 \BOX
\end{remark}

\begin{example}[affine $\lie{gl}_m$ case]
\label{example:gl}
\normalfont
 Assume $m\in\Z_{\geqq2}$ and set $I=\{0,1,\ldots,m-1\}$.
 Let $Y$ be the free $\Z$-module generated 
 by $\{\eps_i\}_{i=1}^m$, $c$, and $d$.
 Let $X$ be the dual lattice of $Y$.
 Define the non-degenerate symmetric bilinear form $(\,,):Y\times Y\to \Z$ 
 by $(\eps_i,\eps_j)=\delta_{ij}$, $(c,d)=1$,
 and $(\eps_i,c)=(\eps_i,d)=(c,c)=(d,d)=0$.
 We can identify $X$ with $Y$ by $(\,,)$.
 We set 
 $\av_i=\alpha_i=\eps_i-\eps_{i-1}$ for $i=1,\ldots,m-1$,
 $\av_0=\alpha_0=c-\eps_1+\eps_m$.
 Define the matrix $A^{(1)}_{m-1}$ by 
 $A^{(1)}_{m-1}=[a_{ij}]=[(\av_i,\alpha_j)]_{i,j\in I}$.
 If $m=2$, then $a_{00}=a_{11}=2$ and $a_{01}=a_{10}=-2$.
 If $m\geqq3$, then \(
  a_{ij} = 2\delta_{ij} 
          -\delta_{i+1,j}-\delta_{j+1,i}
          -\delta_{i0}\delta_{j,m-1}-\delta_{j0}\delta_{i,m-1}
 \).
 Thus $A^{(1)}_{m-1}$ is a symmetric GCM and 
 $(Y,\ab X,\ab\bra\,,\ket,\ab\{\av_i\}_{i\in I},\ab\{\alpha_i\}_{i\in I})$ 
 is a $Y$-regular and $X$-regular root datum of type $A^{(1)}_{m-1}$.
 The Kac-Moody algebra associated to the root datum
 can be identified with the affine Lie algebra 
 $\widehat{\lie{gl}}_m=\lie{gl}_m(\C[t,t^{-1}])\oplus\C c\oplus\C td/dt$ by
 \begin{align*}
  &
  E_0 = t E_{m1}, \ F_0 = t^{-1} E_{1m}, 
  \\ &
  E_i = E_{i,i+1}, \ F_i = E_{i+1,i} 
  \quad \text{for $i=1,\ldots,m-1$},
  \\ &
  \eps_i = E_{ii} \quad\text{for $i=1,\ldots,m$},
  \quad d = td/dt,
 \end{align*}
 where $E_{ij}$ ($i,j=1,\ldots,m$) are unit matrices.
 We set all $d_i=1$.
 The quantized universal enveloping algebra associated to the root datum
 is called {\em the quantized universal enveloping algebra of $\lie{gl}_m$}
 and denoted by $U_q(\widehat{\lie{gl}}_m)$.
 \BOX
\end{example}

\begin{example}[affine $\lie{sl}_m$ and $\lie{psl}_m$ cases]
\label{example:sl-psl}
\normalfont
 Assume $m\in\Z_{\geqq2}$ and set $I=\{0,1,\ldots,m-1\}$.
 Let $A^{(1)}_{m-1}=[a_{ij}]_{i,j\in I}$ be the symmetric GCM given above.
 We set all $d_i=1$.
 Let $Y$ be the free $\Z$-module generated 
 by $\{\av_i\}_{i=1}^{m-1}$, $c$, and $d$.
 Set $\av_0=c-\av_1\cdots-\av_{m-1}$.
 Let $X$ be the dual lattice of $Y$.
 Define the non-degenerate symmetric bilinear form $(\,,):Y\times Y\to \Z$ 
 by $(\av_i,\av_j)=a_{ij}$, $(c,d)=1$,
 and $(\av_i,c)=(\av_i,d)=(c,c)=(d,d)=0$ for $i=1.\ldots,m-1$.
 Identifying $X$ with a sublattice of $\h=\C\o_\Z Y$ by $(\,,)$,
 we have $Y\subsetneqq X$.
 Define $\alpha_j\in X$ for $j\in I$ by
 $\bra\av_i,\alpha_j\ket=a_{ij}$,
 $\bra c,\alpha_j\ket=0$,
 $\bra d,\alpha_j\ket=\delta_{j0}$
 for $i,j\in I$.
 Then 
 $(Y,\ab X,\ab\bra\,,\ket,\ab\{\av_i\}_{i\in I},\ab\{\alpha_i\}_{i\in I})$ 
 is a $Y$-regular and $X$-regular root datum of type $A^{(1)}_{m-1}$.
 The Kac-Moody algebra associated to the root datum
 can be identified with the affine Lie algebra 
 $\widehat{\lie{sl}}_m=\lie{sl}_m(\C[t,t^{-1}])\oplus\C c\oplus\C td/dt$ by
 \begin{align*}
  &
  E_0 = t E_{m1}, \ F_0 = t^{-1} E_{1m}, 
  \\ &
  E_i = E_{i,i+1}, \ F_i = E_{i+1,i} 
  \quad \text{for $i=1,\ldots,m-1$},
  \\ &
  \av_i = E_{ii}-E_{i+1,i+1}
  \quad \text{for $i=1,\ldots,m-1$}, 
  \\ &
  \av_0 = c - E_{11} + E_{mm},
  \quad d = td/dt.
 \end{align*}
 The quantized universal enveloping algebra associated to the root datum
 is called {\em the quantized universal enveloping algebra of $\lie{sl}_m$}
 and denoted by $U_q(\widehat{\lie{sl}}_m)$.
 Associating to the dual root datum 
 $(X,\ab Y,\ab\bra\,,\ket,\ab\{\alpha_i\}_{i\in I},\ab\{\av_i\}_{i\in I})$, 
 we define the quantized universal enveloping algebra 
 $U_q(\widehat{\lie{psl}}_m)$.
 \BOX
\end{example}


\subsection{Adjoint action}
\label{sec:ad}

For an arbitrary Hopf algebra $H$,
the adjoint action $\ad:H\to\End H$ is defined 
by $\ad(x)y=\sum_{(x)}x_{(1)}yS(x_{(2)})$ for $x,y\in H$,
where $\Delta(x)=\sum_{(x)}x_{(1)}\otimes x_{(2)}$ (the Swedler notation).

In a quantized universal enveloping algebra $U_q$, we have
\begin{align*}
 &
 \ad(E_i)x = E_i x K_i^{-1} - x E_i K_i^{-1}
 \quad\text{for $i\in I$, $x\in U_q$},
 \\ &
 \ad(F_i)x = F_i x  - K_i^{-1} x K_i F_i
 \quad\text{for $i\in I$, $x\in U_q$},
 \\ &
 \ad(q^\lambda)x = q^\lambda x q^{-\lambda}
 \quad\text{for $\lambda\in d^{-1}Y$, $x\in U_q$}.
\end{align*}
Setting $x=F_i\o1$ and $y=K_i^{-1}\otimes F$, we have $yx=q_i^2xy.$
Using the $q$-binomial theorem, we obtain \(
  \Delta(F_i)^n 
  \ab = \ab
  {\displaystyle\sum_{k=0}^n} \ab
  q_i^{k(n-k)} \ab \qbinom{n}{k}_{q_i} \ab 
  F_i^{n-k} \ab K_i^{-k}\ab\otimes\ab K_i^k
\) and hence \(
  (1\otimes S)(\Delta(F_i^n))
  \ab = \ab
  {\displaystyle \sum_{k=0}^n} \ab
  (-1)^k \ab q_i^{k(n-1)} \ab \qbinom{n}{k}_{q_i} \ab 
  F_i^{n-k} \ab K_i^{-k} \ab \otimes \ab K_i^k \ab F_i^k
\).  We conclude that
\begin{equation*}
 \ad(F_i)^n x
 = \ad(F_i^n) x
 = \sum_{k=0}^n (-1)^k q_i^{k(n-1)} \qbinom{n}{k}_{q_i}
   F_i^{n-k} K_i^{-k} x K_i^k F_i^k.
\end{equation*}
In particular, we have
\begin{equation}
 \ad(F_i)^n F_j
 = \sum_{k=0}^n
 (-1)^k q_i^{k(n-1+a_{ij})} \qbinom{n}{k}_{q_i}
  F_i^{n-k} F_j F_i^k.
 \label{eq:ad(F_i)^n}
\end{equation}
Hence the $q$-Serre relations for $F_i$ ($i\in I$) are
rewritten as $\ad(F_i)^{1-a_{ij}}F_j=0$ for $i\ne j$.
Similar results hold for $E_i$ ($i\in I$).
The following lemma is equivalent to 
Formula (14) of \cite{iohara-malikov-GKconj}.

\begin{lemma}[\cite{iohara-malikov-GKconj}]
\label{lemma:F_i^nF_j}
 Assume that $x\in U_q$ and $K_i^{-1} x K_i = q_i^a x$.  
 Then
 \begin{equation*}
  F_i^n x 
  = 
  \sum_{k=0}^{\infty} 
  q_i^{(k+a)(n-k)}\qbinom{n}{k}_{q_i}
  \ad(F_i)^k(x)F_i^{n-k}
  \quad\text{for $n\in\Z_{\geqq0}$},
 \end{equation*}
 where the left-hand side is a finite sum with respect to $k=0,1,\ldots,n$.
 In particular, if $i\ne j$ in $I$, then
 \begin{equation*}
  F_i^n F_j 
  = 
  \sum_{k=0}^{-a_{ij}} 
  q_i^{(k+a_{ij})(n-k)}\qbinom{n}{k}_{q_i}
  \ad(F_i)^k(F_j)F_i^{n-k}
  \quad\text{for $n\in\Z_{\geqq0}$}.
  \BOX
 \end{equation*}
\end{lemma}
The first formula of this lemma is proved by induction on $n$.
The case for $n=1$ leads to the cases for any $n\in\Z_{>0}$.
The second immediately follows from the $q$-Serre relations.

The following example can be found as 
Formula (24) of \cite{iohara-malikov-GKconj}.
\begin{example}[\cite{iohara-malikov-GKconj}]
\label{example:F_i^nF_j}
\normalfont
 If $a_{ij}=-1$, then
 \begin{align*}
  F_i^n F_j 
  &= q_i^{-n}F_jF_i^n+[n]_{q_i}\ad(F_i)(F_j)F_i^n-1
  \\ &
  = [1-n]_{q_i}F_jF_i^n + [n]_{q_i}F_iF_jF_i^{n-1}.
  \BOX
 \end{align*}
\end{example}


\subsection{Verma relations}
\label{sec:Verma-relations}

Let $A=[a_{ij}]_{i,j\in I}$ be 
a symmetrizable GCM symmetrized by $\{d_i\}_{i\in I}$
and $(Y,\ab X,\ab\bra\,,\ket,\ab\{\av_i\}_{i\in I},\ab\{\alpha_i\}_{i\in I})$
a $Y$-regular and $X$-regular root datum of type $A$.

Denote by $W$ the Weyl group of type $A$
and by $s_i$ ($i\in I$) its generators.
Denote by $U_q^-$ the lower part of 
the quantized universal enveloping algebra associated to the root datum
and by $F_i$ ($i\in I$) its Chevalley generators.
Let $\lambda\in X^+$.

Assume that $s_{i_1}s_{i_2}\cdots s_{i_n}$ is a reduced expression in $W$.
We set $k_p\in\Z$ for $p=1,2,\ldots,n$ by
\begin{equation*}
 k_p = \bra s_{i_n}s_{i_{n-1}}\cdots s_{i_{p+1}}(\av_{i_p}),\lambda\ket.
\end{equation*}
For examples, $k_n=\bra\av_{i_n},\lambda\ket$, 
$k_{n-1}=\bra s_{i_n}(\av_{i_{n-1}}),\lambda\ket$, 
$k_{n-2}=\bra s_{i_n} s_{i_{n-1}}\ab(\av_{i_{n-2}}),\ab\lambda\ket$, 
and so on.
Since $s_{i_n}\cdots s_{i_{p+1}}s_{i_p}$ is also a reduced expression,
$s_{i_n}\cdots s_{i_{p+1}}\ab(\av_{i_p})\ab\in\sum_{i\in I}\Z_{\geqq0}\av_i$.
Therefore $k_p\in\Z_{\geqq0}$ for $p=1,2,\ldots,n$.

Assume that $s_{j_1}s_{j_2}\cdots s_{j_n}$ is 
another reduced expression with 
$s_{i_1}\ab s_{i_2}\ab\cdots\ab s_{i_n}\ab
=\ab s_{j_1}\ab s_{j_2}\ab\cdots\ab s_{j_n}$.
We similarly set $l_p\in\Z$ for $p=1,2,\ldots,n$ by
\begin{equation*}
 l_p = \bra s_{j_n}s_{j_{n-1}}\cdots s_{j_{p+1}}(\av_{j_p}),\lambda\ket.
\end{equation*} 
Then we have the following identity in $U_q^-$:
\begin{equation}
   F_{i_1}^{k_1}F_{i_2}^{k_2}\cdots F_{i_n}^{k_n}
 = F_{j_1}^{l_1}F_{j_2}^{l_2}\cdots F_{j_n}^{l_n}.
 \label{eq:Verma1}
\end{equation}
Furthermore the sequence $((i_1,\ab k_1),\ab(i_2,\ab k_2),\ab\ldots,\ab(i_n,\ab k_n))$ 
is equal to the sequence $((j_1,\ab l_1),\ab(j_2,\ab l_2),\ab\ldots,\ab(j_n,\ab l_n))$ 
up to permutation of order.
In order to prove these results, 
it is sufficient to show them for each pair of reduced expressions 
in \exampleref{example:braid-relations}.
For the proof,
see Section 39.3 of \cite{lusztig} and 
Lemma 2 of \cite{etingof-varchenko-dynamical-Weyl-1}.
These results are called {\em Verma relations}.

\begin{example}[Verma relations]
\label{example:Verma-relations}
\normalfont
 Assume that $i,j\in I$ and $i\ne j$.
 Let $k$ and $l$ be arbitrary non-negative integers. 
 Then \exampleref{example:braid-relations} leads to 
 the following formulae of the Chevalley generators of $U_q^-$:
 \begin{enumerate}
  \item $F_i^{k}F_i^{l} = F_j^{l}F_i^{k}$
	if $(a_{ij},a_{ji})=(0,0)$;
  \item $F_i^{l}F_j^{k+l}F_i^{k} = F_j^{k}F_i^{k+l}F_j^{l}$
	if $(a_{ij},a_{ji})=(-1,-1)$;
  \item $F_i^{k}F_j^{2k+l}F_i^{k+l}F_j^{l} = F_j^{l}F_i^{k+l}F_j^{2k+l}F_i^{k}$
	if $(a_{ij},a_{ji})=(-1,-2)$;
  \item $F_i^{k}F_j^{3k+l}F_i^{2k+l}F_j^{3k+2l}F_i^{k+l}F_j^{l} 
       = F_j^{l}F_i^{k+l}F_j^{3k+2l}F_i^{2k+l}F_j^{3k+l}F_i^{k}$
	if $(a_{ij},a_{ji})=(-1,-3)$.
 \end{enumerate}
 These formulae shall be used in the construction of 
 the quantized birational Weyl group actions.
 \BOX
\end{example}

\begin{remark}
\label{remark:Verma-relations}
\normalfont
 For $\beta_1,\ldots,\beta_n\in\sum_{i\in I}\Z_{\geqq0}\av_i$,
 denote by $F_{i_1}^{\beta_1}\cdots F_{i_n}^{\beta_n}$ 
 the mapping from $X^+$ to $U_q^-$ sending $\lambda$ to
 $F_{i_1}^{\bra\beta_1,\lambda\ket}\cdots F_{i_n}^{\bra\beta_n,\lambda\ket}$.
 We introduce the formal symbols $\ts_i$ ($i\in I$)
 satisfying the braid relations and 
 $\ts_i^{-1} F_j^{\beta_j} \ts_i = F_j^{s_i(\beta_j)}$.
 Then the Verma identity \eqref{eq:Verma1} can be formally rewritten 
 in the following form:
 \begin{equation*}
  \ts_{i_1}F_{i_1}^{\av_{i_1}}
  \ts_{i_2}F_{i_2}^{\av_{i_2}}
  \cdots 
  \ts_{i_n}F_{i_n}^{\av_{i_n}}
  = 
  \ts_{j_1}F_{j_1}^{\av_{j_1}}
  \ts_{j_2}F_{j_2}^{\av_{j_2}}
  \cdots 
  \ts_{j_n}F_{j_n}^{\av_{j_n}}.
 \end{equation*}
 This means that
 $\ts_i F_i^{\av_i}$ ($i\in I$) formally satisfy the braid relations.
 Moreover, if we have $\ts_i^2=1$ and $F_i^{-\av_i}F_i^{\av_i}=1$, 
 then we obtain, at least formally, 
 \begin{equation*}
    \ts_i F_i^{\av_i} \ts_i F_i^{\av_i}
  = \ts_i^2 F_i^{-\av_i} F_i^{\av_i}
  = 1.
 \end{equation*}
 This means that
 $\ts_i F_i^{\av_i}$ ($i\in I$) formally satisfy
 the defining relations of the Weyl group.
 If we can justify the above heuristic consideration, 
 then we can construct the braid or Weyl group representations.
 \BOX
\end{remark}


\section{Localizations of non-commutative rings}
\label{sec:localizations}

In this section, we shall summarize results on
localizations of non-commutative rings
necessary to quantize birational actions.
Most of the proofs omitted below can be found, for examples, 
in Chapter 10 of \cite{goodearl-warfield-ncalg}
and Chapter 2 of \cite{mcconnell-robson-ncalg}.


\subsection{Localization at an Ore subset}
\label{sec:Ore-subsets}

Let $A$ be a (possibly non-commutative) ring.
$A$ is called an {\em integral domain} (or a {\em domain} for short)
if $A\ne 0$ and the products of non-zero elements of $A$ are always non-zero.
A proper two-sided ideal $I$ of $A$ is called {\em completely prime}
if $A/I$ is an integral domain.
We say that $A$ is {\em left Noetherian}
if there is no infinite properly ascending chain of left ideals of $A$.
A {\em right Noetherian} ring is similarly defined.

A subset $S$ of $A$ is called {\em multiplicative} 
if $S$ contains $1$ and is closed with respect to multiplication.
Let $A$ be an integral domain and $S$ its multiplicative subset.
We say that $S$ satisfies 
{\em the left} (resp.\ {\em right}) {\em Ore condition}
if $Sa\cap As\ne\emptyset$ (resp.\ $aS\cap sA\ne\emptyset$) 
for any $a\in A$ and $s\in S$.
A multiplicative subset satisfying 
the left (resp.\ right) Ore condition
is called a {\em left} (resp.\ {\em right}) {\em Ore subset} for short.
A left and right Ore subset is simply called an {\em Ore subset}.

Assume that $S$ is a left Ore subset of $A$.
Then we can define the ring $S^{-1}A$ as follows.
As a set, $S^{-1}A$ is defined to be 
the quotient set $(S\times A)/{\sim}$, 
where the equivalence relation $\sim$ is defined by
$(s,a)\sim(s',a')$ $\Leftrightarrow$ 
there exists $u,u'\in A$ such that $us=u's'\in S$ and $ua=u'a'$.
Denote by $s\backslash a$ 
the element of $S^{-1}A$ represented by $(s,a)\in S\times A$.
We can define the ring structure of $S^{-1}A$ by
\begin{alignat*}{2}
 &
 (s\backslash a)(s'\backslash a') = (s''s)\backslash(a''a'),
 & \quad & s''a = a''s',\ a''\in A,\ s''\in S;
 \\ &
 s\backslash a + s'\backslash a' = (us)\backslash(ua+u'a'),
 & \quad & us=u's',\ u'\in A,\ u\in S.
\end{alignat*}
Identifying $a\in A$ with $1\backslash a\in S^{-1}A$, 
we can embed $A$ into $S^{-1}A$.
Then the ring $S^{-1}A$ contains $A$ as a subring
and satisfies that
any element of $S$ is invertible in $S^{-1}A$ and
$S^{-1}A=\{\,s^{-1}a=s\backslash a\mid s\in S,\, a\in A\,\}$.
Furthermore $S^{-1}A$ has the following universality:
for any ring $B$ and any ring homomorphism $f:A\to B$ with
the property that $f(s)$ is invertible in $B$ for any $s\in S$,
there exists a unique ring homomorphism $\phi:S^{-1}A\to B$ 
with $\phi|_A=f$.
In particular $S^{-1}A$ is uniquely, up to isomorphism, 
determined by $A$ and $S$.
We call $S^{-1}A$ the {\em left localization} of $A$ at $S$.
If $S$ is an Ore subset, then the left localization at $S$ 
can be identified with the right one, 
namely $S^{-1}A=\{\,as^{-1}\mid s\in S,\, a\in A\,\}$.

\begin{lemma}
\label{lemma:Ore-subset}
 Let $A$ be an integral domain 
 generated by $\{a_j\}_{j\in J}$ over a field $\F$
 and $S$ its multiplicative subset generated by $\{s_i\}_{i\in I}$.
 Then we have the following results:
 \begin{enumerate}
  \item \label{enumi:Ore-subset-1}
   If $Sa\cap As_i\ne\emptyset$ for any $a\in A$ and $i\in I$,
   then $S$ is a left Ore subset of $A$.
  \item \label{enumi:Ore-subset-2}
   Assume that for any $i\in I$, $j\in J$, and $n\in\Z_{>0}$,
   there exists $N\in\Z_{>0}$ with $s_i^N a_j\in A s_i^n$.
   Then for any $i\in I$, $a\in A$, and $n\in\Z_{>0}$,
   there exists $N\in\Z_{>0}$ with $s_i^N a  \in A s_i^n$.
   Therefore $S$ is a left Ore subset of $A$.
 \end{enumerate}
\end{lemma}

\begin{proof}
 \enumiref{enumi:Ore-subset-1}
 Take $i_1,\ldots,i_n\in I$.
 By induction on $n$, 
 let us show that $Sa\cap Af_{i_n}\cdots f_{i_1}\ne\emptyset$ 
 for any $a\in A$.
 The case of $n=1$ is just the the assumption.
 Assume that it holds for $n-1$.
 Then there exist $t\in S$ and $b\in A$ with $ta=bf_{i_{n-1}}\cdots f_{i_1}$.
 By the case of $n=1$,
 there exist $c\in A$ and $u\in S$ with $ub=cf_{i_n}$.
 Then $ut\in S$ and $uta=cf_{i_n}\cdots f_{i_1}$.

 \enumiref{enumi:Ore-subset-2}
 Fix any $i\in I$.
 Let $\widetilde{A}$ be the subset of $A$ 
 consisting of the elements $a\in A$ such that
 for any $n\in\Z_{>0}$ there exists $N\in\Z_{>0}$ with $s_i^N a\in A s_i^n$.
 It is sufficient 
 for the proof of the first statement
 to show that $\widetilde{A}$ is a subalgebra of $A$.
 Take any $a,b\in \widetilde{A}$.
 For any $n\in\Z_{>0}$,
 there exists $M,N\in\Z_{>0}$ such that 
 $s_i^M a\in A s_i^n$ and $s_i^N b\in A s_i^n$.
 Then $s_i^{M+N}(a+b)\in As_i^n$ and hence $a+b\in \widetilde{A}$.
 There exists $L\in\Z_{>0}$ such that $s_i^La\in As_i^N$.
 Then $s_i^Lab\in A s_i^N b\subset A s_i^n$ and hence $ab\in \widetilde{A}$.
 We have shown that $\widetilde{A}$ is a subalgebra of $A$.
 The second statement follows from \enumiref{enumi:Ore-subset-1}.
\end{proof}

\begin{example}[the inverse of $F_i$]
\label{example:Ore-subset-of-U_q^-}
\normalfont
 Consider a quantized universal enveloping algebra $U_q$
 and its lower part $U_q^-$.
 Let $J$ be any subset of $I$
 and $S_J$ the multiplicative subset generated by $\{F_j\}_{j\in J}$.
 Using Formula \eqref{eq:F^nE^m},
 $F_iq_i^\lambda=q_i^{\lambda+\bra\lambda,\alpha_i\ket}F_i$ 
 ($\lambda\in d^{-1}Y$), 
 and the $q$-Serre relations of $\{F_i\}_{i\in I}$, 
 we can find that
 both $(U_q^-,S_J)$ and $(U_q,S_J)$ satisfy the assumption of 
 \lemmaref{lemma:Ore-subset} \enumiref{enumi:Ore-subset-2}.
 Therefore $S_J$ is a left Ore subset of $U_q^-$ and $U_q$.
 The anti-algebra involution
 given by $E_i\mapsto E_i$, $F_i\mapsto F_i$ ($i\in I$)
 and $q^\lambda\mapsto q^{-\lambda}$ ($\lambda\in d^{-1}Y$)
 proves that $S_J$ is also a right Ore subset of $U_q^-$ and $U_q$.
 By the universality of $S_J^{-1}U_q^-$,
 we can regard $S_J^{-1}U_q^-$ as a subalgebra of $S_J^{-1}U_q$.
 For $J=\{i_1,\ldots,i_r\}$, 
 we denote $S_J^{-1}U_q$ by $U_q[F_{i_1}^{-1},\ldots,F_{i_r}^{-1}]$
 and $S_J^{-1}U_q^-$ by $U_q^-[F_{i_1}^{-1},\ldots,F_{i_r}^{-1}]$.
 \BOX
\end{example}

Using the inverse of $F_i$, we can state the following generalization of 
\lemmaref{lemma:F_i^nF_j}
for negative integral powers of $F_i$.

\begin{lemma}
\label{lemma:F_i^nF_j-all}
 Assume that $x\in U_q$, $K_i^{-1} x K_i = q_i^a x$,
 and $\ad(F_i)^kx=0$ for sufficiently large $k$.  
 Then we have the following formula in $U_q[F_i^{-1}]$:
 \begin{equation*}
  F_i^n x 
  = 
  \sum_{k=0}^{\infty} 
  q_i^{(k+a)(n-k)}\qbinom{n}{k}_{q_i}
  \ad(F_i)^k(x)F_i^{n-k}
  \quad\text{for $n\in\Z$},
 \end{equation*}
 where the left-hand side is a finite sum.
 In particular, if $i\ne j$ in $I$, then
 \begin{equation*}
  F_i^n F_j 
  = 
  \sum_{k=0}^{-a_{ij}} 
  q_i^{(k+a_{ij})(n-k)}\qbinom{n}{k}_{q_i}
  \ad(F_i)^k(F_j)F_i^{n-k}
  \quad\text{for $n\in\Z$}.
  \BOX
 \end{equation*}
\end{lemma}

\begin{proof}
 The second formula immediately follows from 
 the first formula and the $q$-Serre relations.
 By \lemmaref{lemma:F_i^nF_j}, we can assume that $n$ is negative.
 By induction on $N\in\Z_{>0}$, we can obtain the following formula:
 \begin{align*}
  F_i^{-1} x 
  &
  = \sum_{k=0}^{N-1} (-1)^k q_i^{-(k+1)(k+a)} \ad(F_i)^k(x)F_i^{-(k+1)}
  \\ &
  + q_i^{-N(N-1+a)} F_i^{-1}\ad(F_i)^N(x)F_i^{-N}.
 \end{align*} 
 Since $\ad(F_i)^k(x)=0$ for sufficiently large $k$
 and $\qbinom{-1}{k}_{q_i}=(-1)^k$, 
 the first formula for $n=-1$ has been proved.
 This leads to the first formulae for all negative $n$ by induction on $-n$.
\end{proof}


\subsection{Ore domains}
\label{sec:Ore-domains}

An integral domain $A$ is called 
a {\em left} (resp.\ {\em right}) {\em Ore domain}
if $Aa\cap Ab\ne 0$ (resp.\ $aA\cap bA\ne 0$) for any non-zero $a,b\in A$.
In other words, an integral domain $A$ is a left (resp.\ right) Ore domain
if and only if $A\setminus\{0\}$ is a left (resp.\ right) Ore subset.
A left and right Ore domain is simply called an {\em Ore domain}.

Assume that $A$ is an Ore domain.
Let $K$ be the localization of $A$ at $A\setminus\{0\}$.
Then $K$ is a skew field and \(
 K = \{\, s^{-1}a \mid a,s\in A,\ s\ne0 \,\}
   = \{\, as^{-1} \mid a,s\in A,\ s\ne0 \,\}
\). We call $K$ the {\em (skew) field of fractions} of $A$
and denote $K$ by $Q(A)$.

\begin{lemma}[2.1.15 of \cite{mcconnell-robson-ncalg}]
\label{lemma:noetherian-domains}
 A left (resp.\ right) Noetherian domain 
 is a left (resp.\ right) Ore domain.
 In particular a left and right Noetherian domain is an Ore domain.
 \BOX
\end{lemma}

\begin{example}
\label{example:noetherian-domains}
\normalfont
 The following are left and right Noetherian domains
 (Chapter 1 of \cite{mcconnell-robson-ncalg}):
 \begin{enumerate}
 \item the skew polynomial ring $R[x;\sigma,\delta]$ 
  associated to 
  a left and right Noetherian domain $R$, 
  an algebra automorphism $\sigma$ of $R$,
  and a $\sigma$-derivation $\delta$ of $R$ 
  ($\delta(ab)=\delta(a)b+\sigma(a)\delta(b)$ for $a,b\in R$),
  defined to be the ring generated by $a\in R$ and $x$
  with defining relations: 
  (a) $R$ is a subring of $R[x;\sigma,\delta]$,
  (b) $xa=\sigma(a)x+\delta(a)$ for $a\in R$;
 \item the skew Laurent polynomial ring $R[x,x^{-1};\sigma]$ 
  associated to 
  a left and right Noetherian domain $R$
  and an algebra automorphism $\sigma$ of $R$,
  defined to be the ring generated by $a\in R$ and $x^{\pm1}$
  with defining relations:
  (a) $R$ is a subring of $R[x,x^{-1};\sigma]$,
  (b) $xa=\sigma(a)x$ for $a\in R$,
  (c) $xx^{-1}=x^{-1}x=1$;
 \item the Weyl algebras over a filed $\F$ of characteristic $0$
  generated by $x_1,\ldots,x_n,\d_1,\ldots,\d_n$
  with defining relations:
  $x_ix_j=x_jx_i$, $\d_i\d_j=\d_j\d_i$, 
  and $\d_i x_j - x_j \d_i = \delta_{ij}$;
 \item the universal enveloping algebra $U(\g)$ of 
  any finite dimensional Lie algebra $\g$ over a field.
 \end{enumerate}
 Let $\F$ be a field and $q_{ij}\in\F^\times$ for $i,j=1,\ldots,n$.
 Assume that $q_{ii}=1$ and $q_{ji}=q_{ij}^{-1}$.
 From the first and  second examples above, we obtain, by induction on $n$, 
 the following examples of left and right Noetherian domains
 respectively:
\begin{enumerate}
 \setcounter{enumi}{4}
 \item the $q$-polynomial ring over $\F$
  defined to be the algebra over $\F$ generated by
  $x_1,\ldots,x_n$ with defining relations
  $x_jx_i=q_{ij}x_ix_j$ for any $i,j$;
 \item the $q$-Laurent polynomial ring over $\F$
  defined to be the algebra over $\F$ generated by
  $x_1^{\pm1},\ldots,x_n^{\pm1}$ with defining relations
  $x_i^{-1}x_i=x_ix_i^{-1}=1$ for any $i$
  and $x_jx_i=q_{ij}x_ix_j$ for any $i,j$.
 \end{enumerate}
 All of these examples are Ore domains.
 \BOX
\end{example}

In the next subsection, we shall deal with
Ore domains which are not always left and right Noetherian.


\subsection{Tempered domains}
\label{sec:slowly-increasing}

Let $A$ be an associative algebra over a field $\F$ and 
$\{F_k A\}_{k=0}^\infty$ a family of $\F$-vector subspaces of $A$.
Set $A_k=0$ for $k\in\Z_{<0}$.
We say that $\{F_k A\}_{k=0}^\infty$ is 
a {\em filtration} of $A$
if $F_0 A\subset F_1 A\subset F_2 A\subset\cdots$,
$\bigcup_{k=0}^\infty F_k A=A$,
$1\in F_0 A$, and $F_k A\,F_l A\subset F_{k+l}A$ for any $k,l$.
Let $\{F_k A\}_{k=0}^\infty$ be a filtration of $A$.
Set $\gr_k A=F_k A/F_{k-1}A$ and $\gr A=\bigoplus_{k=0}^\infty\gr_k A$.
Then $\gr A$ has a natural graded algebra structure.
If $\gr A$ is an integral domain 
(resp.\ left Noetherian, right Noetherian), 
then $A$ is so.

\begin{definition}[tempered domain]
\label{definition:tempered-domain}
\normalfont
 An associative algebra $A$ over a field $\F$
 has a {\em slowly increasing filtration} 
 if there exists a filtration $\{F_k A\}_{k=0}^\infty$ of $A$ 
 such that $\limsup_k\ab(\dim_\F F_k A)^{1/k}\leqq 1$.
 This is equivalent to the condition that 
 $\dim_\F\gr_k A<\infty$ for all $k$ and 
 $\limsup_k\ab(\dim_\F\gr_k A)^{1/k}\leqq 1$.
 See \remarkref{remark:convergence-radius} \enumiref{enumi:tempered-1} below.
 An associative algebra with slowly increasing filtration
 is called a {\em tempered algebra} for short.
 In addition, if $A$ is an integral domain, 
 then $A$ is called a {\em tempered domain}.
 \BOX
\end{definition}

From the definition we can immediately obtain the following result.

\begin{lemma}
\label{lemma:subquot-tempered}
 Assume that $A$ is a tempered domain over a field. 
 Then subalgebras of $A$ and 
 quotient integral domains of $A$
 are also tempered domains.
 \BOX
\end{lemma}

\begin{remark}
\label{remark:convergence-radius}
\normalfont
 Let $\{a_k\}_{k=0}^\infty$ be a sequence of complex numbers
 and $\rho$ the convergence radius of 
 the power series $\sum_{k=0}^\infty a_k z^k$.
 The Cauchy-Hadamard theorem says 
 $\limsup_k\ab|a_k|^{1/k}=\rho^{-1}$.
 The absolute convergence of $\sum_{k=0}^\infty a_k$ is equivalent to
 the condition that, for sufficiently large $k_0$,
 the infinite product $\prod_{k=k_0}^\infty(1+a_k)$ 
 is absolutely convergent to a non-zero complex number.
 Therefore the following conditions are mutually equivalent:
 \begin{itemize}
  \item[\upshape(a)] 
   $\limsup_k\ab|a_k|^{1/k}\leqq 1$.
  \item[\upshape(b)] 
   There exists a holomorphic function in $|z|<1$
   such that its Maclaurin expansion is equal to $\sum_{k=0}^\infty a_k z^k$.
  \item[\upshape(c)] 
   For sufficiently large $k_0$,
   the infinite product $\prod_{k=k_0}^\infty(1+a_kz^k)$ 
   is absolutely and uniformly convergent in wide sense 
   to a non-vanishing holomorphic function in $|z|<1$.
 \end{itemize}
 These observations lead to the following results:
 \begin{enumerate}
  \item \label{enumi:tempered-1}
   $\limsup_k\ab|a_k|^{1/k}\leqq 1$
   implies $\limsup_k\ab|a_0+a_1+\cdots+a_k|^{1/k}\leqq 1$.
  \item \label{enumi:tempered-2}
   Fix $N\in\Z$.
   Then $\limsup_k\ab|a_k|^{1/k}\leqq 1$
   is equivalent to $\limsup_k\ab|a_k|^{1/(k+N)}\leqq 1$.
  \item \label{enumi:tempered-3}
   Assume that $|a_k|\leqq 1$ for all $k$.
   For any positive integer $N$, 
   define the sequence $\{b^{(N)}_k\}_{k=0}^\infty$
   of complex numbers by the Maclaurin expansion \(
    \left(\prod_{k=0}^\infty(1-a_k z^k)\right)^{-N}
    = \sum_{k=0}^\infty b^{(N)}_k z^k
   \) in $|z|<1$. Then $\limsup_k\ab|b_k|^{1/k}\leqq 1$.
  \item \label{enumi:tempered-4}
   Set $c_k=\sum_{i=0}^k a_ib_{k-i}$. Then
   $\limsup_k\ab|a_k|^{1/k}\leqq 1$ and $\limsup_k\ab|b_k|^{1/k}\leqq 1$
   implies $\limsup_k\ab|c_k|^{1/k}\leqq 1$.
   \BOX
 \end{enumerate}
\end{remark}

Rocha-Caridi and Wallach found the following criterion of Ore domains
(Lemma 1.2 of \cite{rochacaridi-wallach-ore}) 
to prove that 
the universal enveloping algebras of
affine (or Euclidean) Lie algebras are Ore domains.

\begin{lemma}[\cite{rochacaridi-wallach-ore}]
\label{lemma:tempered}
 A tempered domain over a field is always an Ore domain.
\end{lemma}

\begin{proof}
 Let $A$ be an integral domain over a field $\F$ and
 $\{F_k A\}_{k=0}^\infty$ an arbitrary increasing filtration of $A$.
 Assume that $A$ is not an Ore domain.
 Then there exist $r\in\Z_{\geqq0}$ and
 non-zero $a,b\in F_r A$ such that $Aa\cap Ab=0$ or $aA\cap bA=0$.
 Assume that $Aa\cap Ab=0$.
 Set $a_k=\dim_\F F_k A$ and 
 take $k_0\in\Z_{\geqq0}$ so that $a_{k_0}\geqq 1$.
 Since $(F_k A)a+(F_k A)b\subset F_{k+r}A$ 
 and $(F_k A)a\cap(F_k A)b=0$,
 we have $a_{k+r}=\dim_\F F_{k+r}A\geqq\dim((F_k A)a\oplus(F_k A)b)=2a_k$
 and hence $a_{k_0+pr}\geqq 2^p a_{k_0}\geqq 2^p$,
 i.e.\ $a_{k_0+pr}^{1/(pr)}\geqq 2^{1/r}>1$ for all $p\in\Z_{\geqq0}$.
 This leads to $\limsup_k\ab a_k^{1/(k-k_0)}>1$ and hence
 \remarkref{remark:convergence-radius} \enumiref{enumi:tempered-2}
 shows $\limsup_k\ab a_k^{1/k}>1$.
 Therefore $A$ is not tempered.
 When $aA\cap bA=0$, similarly $A$ is not.
 We complete the proof of the lemma.
\end{proof}

\begin{example}
\label{example:q-Laurent}
 Let $A$ be a $q$-Laurent polynomial ring over a field.
 (See \exampleref{example:noetherian-domains}.)
 Then $A$ is a tempered domain.
 Therefore its subalgebras and quotient integral domains
 are also tempered domains.
 Furthermore all of these are Ore domains.
 \BOX
\end{example}

\remarkref{remark:convergence-radius} \enumiref{enumi:tempered-4}
immediately lead to the following result.

\begin{lemma}
\label{lemma:tensor-tempered}
 For any tempered algebras $A$ and $B$ over a field $\F$, 
 the tensor product algebra $A\otimes B$ over $\F$ is also tempered.
 Therefore,  
 if $A$ is a tempered domain
 and $B$ is a $q$-polynomial or $q$-Laurent polynomial ring,
 then $A\otimes B$ is also a tempered domain.
 \BOX
\end{lemma}

\begin{theorem}
\label{theorem:U(n),Uq-}
 The following algebras are tempered domains:
 \begin{enumerate}
 \item 
  the universal enveloping algebras $U(\n_\pm)$
  of the upper and lower parts $\n_\pm$
  of a Kac-Moody algebra of finite or affine type
  (Theorem 1.10 of \cite{rochacaridi-wallach-ore}),
 \item 
  the upper and lower parts $U_q^\pm$ 
  of a quantized universal enveloping algebra of finite or affine type.
 \end{enumerate}
 Therefore these are Ore domains.
\end{theorem}

\begin{proof}
 Denote $U(\n_\pm)$ or $U_q^\pm$ by $A$.
 Let us define a filtration of $A$ 
 by using principal gradations.
 Set $F_k A=\bigoplus_{i=0}^k U(\n_\pm)_{\pm i}$ if $A=U(\n_\pm)$
 and $F_k A=\bigoplus_{i=0}^k (U_q^\pm)_{\pm i}$ if $A=U_q^\pm$.
 Then $\{F_k A\}_{k=0}^\infty$ is a filtration of $A$.
 Because of $U(\n_-)$ and $U_q^-$ are of finite or affine type, 
 from \lemmaref{lemma:dimU(n)k} and \lemmaref{lemma:dimUq(n)k}
 together with 
 \remarkref{remark:convergence-radius} \enumiref{enumi:tempered-3},
 we obtain that $\{F_k A\}_{k=0}^\infty$ is slowly increasing.
 This means that $A$ is a tempered domain.
 Therefore \lemmaref{lemma:tempered} completes the proof.
\end{proof}

In \cite{berman-cox-toroidal}, Berman and Cox showed that 
the universal enveloping algebras 
of Kac-Moody Lie algebras of affine type, 
as well as those of toroidal Lie algebras, 
are tempered domains.
Using an associative algebra version of 
Lemma 1.4 (a) in \cite{berman-cox-toroidal},
we shall show that
quantized universal enveloping algebras of affine type
are also tempered domains.

\begin{lemma}
\label{lemma:tempered-triangular}
 Let $A$ be a domain over a field 
 and $A^\pm$ and $A^0$ its subalgebras.
 Assume that $A^\pm$ and $A^0$ have a slowly increasing filtration
 denoted by 
 $\{F_i A^\pm\}_{i=0}^\infty$ and $\{F_i A^0\}_{i=0}^\infty$
 respectively.
 Assume that these satisfy the following conditions:
 \begin{itemize}
  \item[\upshape(a)] 
   The multiplication gives an isomorphism
   $A^-\otimes A^0\otimes A^+\isomto A$ of vector spaces.
  \item[\upshape(b)] 
   $F_j A^0 F_l A^- = F_l A^- F_j A^0$ and 
   $F_k A^+ F_m A^0 = F_m A^0 F_k A^+$ for any $j,k,l,m$.
  \item[\upshape(c)] 
   \( 
     F_k A^+ F_l A^- \subset 
    \sum_{p=0}^{\max\{k,l\}} F_{l-p}A^- F_p A^0 F_{k^p}A^+
   \) for any $k,l$
 \end{itemize}
 Then $A$ is also a tempered domain and hence an Ore domain.
\end{lemma}

\begin{proof}
 Using the above conditions, 
 we can define a filtration $\{F_l A\}_{l=0}^\infty$ of $A$ by
 \(
  F_l A = \sum_{i+j+k=l} F_i A^- F_j A^0 F_k A^+
 \).
 From \remarkref{remark:convergence-radius} \enumiref{enumi:tempered-4}
 we obtain that $\{F_l A\}_{l=0}^\infty$ is slowly increasing
 and hence $A$ is a tempered domain.
 Therefore \lemmaref{lemma:tempered} completes the proof.
\end{proof}

\begin{theorem}
\label{theorem:affine-tempered}
 The following algebras are tempered domains:
 \begin{enumerate}
 \item 
  the universal enveloping algebra $U(\g)$
  of a Kac-Moody algebra
  of finite or affine type
  (part of Proposition 1.7 of \cite{berman-cox-toroidal}),
 \item 
  a quantized universal enveloping algebra $U_q$
  of finite or affine type.
 \end{enumerate}
 Therefore these are Ore domains.
\end{theorem}

\begin{proof}
 Since $U(\g)$ and $U_q$ are integral domains,
 it is sufficient for the proof to construct
 slowly increasing filtrations of $U(\g)$ and $U_q$.
 
 First we assume that $A=U(\g)$, 
 a Kac-Moody algebra of finite or affine type.
 Set $A^\pm=U(\n_\pm)$ and $A^0=U(\h)$.
 Using the principal gradations of $U(\n_\pm)$, 
 we can define the increasing filtrations $\{F_k A^\pm\}_{k=0}^\infty$ 
 of $A^\pm$ by $F_k A^\pm=\bigoplus_{i=0}^k U(\n_\pm)_{\pm k}$.
 Then we can find from the proof of \theoremref{theorem:U(n),Uq-}
 that $\{F_k A^\pm\}_{k=0}^\infty$ are slowly increasing.
 Let $\{y_1,\ldots,y_M\}$ a basis of $\h$.
 Define the degree by $\deg y_i=1$ for $i=1,\ldots,M$.
 Let $F_k A^0$ be the subspace of 
 $A^0=U(\h)=\C[y_1,\ldots,y_M]$ 
 spanned by the elements of degree $\leqq k$.
 Then $\{F_k A^0\}_{k=0}^\infty$ is 
 a slowly increasing filtration of $A^0$.
 These satisfies the conditions (a), (b), and (c) of 
 \lemmaref{lemma:tempered-triangular}.
 Therefore $A=U(\g)$ is a tempered domain.

 Second we assume that $A=U_q$, the quantized universal enveloping algebra
 associated to a root datum 
 $(Y,\ab Z,\ab\bra\,,\ket,\ab\{\av_i\}_{i\in I},\ab\{\alpha_i\}_{i\in I})$ 
 of finite or affine type.
 Using the principal gradations of $U_q^\pm$, 
 we can define the filtrations $\{F_k A^\pm\}_{k=0}^\infty$ 
 of $A^\pm$ by $F_k A^\pm=\bigoplus_{i=0}^k(U_q^\pm)_{\pm k}$.
 Then we can find from the proof of \theoremref{theorem:U(n),Uq-}
 that $\{F_k A^\pm\}_{k=0}^\infty$ are slowly increasing.
 Let $\{y_1,\ldots,y_M\}$ be a $\Z$-free basis of $d^{-1}Y$.
 Define the degree by $\deg q^{\pm y_i}=1$ for $i=1,\ldots,M$.
 Let $F_k A^0$ be the subspace of 
 $A^0=U_q^0=\F[q^{\pm y_1}\ldots,q^{\pm y_M}]$ 
 consisting of the elements of degree $\leqq k$.
 Then $\{F_k A^0\}_{k=0}^\infty$ is 
 a slowly increasing filtration of $A^0$.
 These satisfies the conditions (a), (b), and (c) of 
 \lemmaref{lemma:tempered-triangular}.
 Therefore $A=U_q$ is a tempered domain.
\end{proof}

From \lemmaref{lemma:subquot-tempered} and the above theorem,
we immediately obtain the following.

\begin{cor}
\label{cor:subquot-affine-tempered}
 Let $A$ be the universal enveloping algebra $U(\g)$ of 
 a Kac-Moody algebra of finite or affine type
 or a quantized universal enveloping algebra $U_q$ of finite or affine type.
 Assume that $B$ is a subalgebra of $A$ and $I$ is a completely prime 
 ideal of $B$.
 Then $B/I$ is a tempered domain and hence an Ore domain.
 \BOX
\end{cor}


\subsection{Truncated $q$-Serre relations}
\label{sec:truncated-q-Serre}

In this subsection, we shall explain a method for constructing 
quotient tempered domains of $U(\n_-)$ and $U_q^-$ for any symmetrizable GCM.

First let us consider the case of $q=1$.

Let $A=\{a_{ij}\}_{i,j\in I}$ be a symmetrizable GCM 
symmetrized by $\{d_i\}_{i\in I}$,
$(Y,\ab X,\ab\bra\,,\ket,\ab\{\av_i\}_{i\in I},\ab\{\alpha_i\}_{i\in I})$ 
a root datum of type $A$,
and $\g$ the Kac-Moody algebra associated to the root datum.
Denote by $\n_-$ (resp.\ $\b_-$) the lower part
(resp.\ the lower Borel subalgebra) of $\g$.
Let $\{y_1,\ldots,y_M\}$ be a basis of $\h=\C\o_\Z Y$.
Assume that if $i\ne j$ and $a_{ij}\ne 0$, 
then $\epsilon_{ij}=\pm1$ and $\epsilon_{ji}=-\epsilon_{ij}$, 
otherwise $\epsilon_{ij}=0$.

Define the algebra $\B$ 
to be the associative algebra over $\C$
generated by $f_i$ ($i\in I$) and $h\in\h$
with defining relations:
\begin{align*}
 &
 \text{$U(\h)=S(\h)$ is a subalgebra of $\B$};
 \\ &
 [h,f_i] = -\bra h,\alpha_i\ket f_i 
 \quad\text{for $i\in I$, $h\in\h$};
 \\ &
 [f_i,f_j] = -\epsilon_{ij}d_ia_{ij}
 \quad\text{for $i,j\in I$}.
\end{align*}
The last relations are sufficient conditions 
of Serre relations for $\{f_i\}_{i\in I}$
and called {\em truncated Serre relations}.
Sending $F_i$ to $f_i$ for each $i\in I$,
we can regard $\B$ as a quotient algebra of $U(\b_-)$.

Define the degree by $\deg f_i=\deg h=1$ for $i\in I$ and $h\in\h$.
Let $F_k\B$ be the subspace of $\B$ 
spanned by the elements of degree $\leqq k$.
Then $\{F_k\B\}_{k=0}^\infty$ is a slowly increasing filtration
of $\B$
and $\gr\B=\bigoplus_{k=0}^\infty F_k\B/F_{k-1}\B$ 
is isomorphic to the commutative polynomial ring
generated by $f_i$ ($i\in I$) and $y_\mu$ ($\mu=1,\ldots,M$).
Therefore $\B$ is a tempered domain.
By \lemmaref{lemma:tensor-tempered}, 
$\B^{\otimes N}$ is also a tempered domain for any positive integer $N$.

We can define the algebra homomorphism $\phi_N:U(\n_-)\to\B^{\otimes N}$
by $\phi_N(F_i)=\sum_{\nu=1}^N f_{i\nu}$,
where $f_{i\nu}=1^{\o(\nu-1)}\otimes F_i\o1^{\o(N-\nu)}$.
Denote the image of $\phi_N$ by $\N_N$.
Then $\N_N$ is also a tempered domain and hence an Ore domain.
Denote $\N_1$ by $\N$ for short.

Second let us consider a $q$-analogue of the above construction.

Let $d$ be the least common denominator of $\{d_i\}_{i\in I}$.
Set the base field $\F$ by $\F=\Q(q^{1/d})$
and $q_i\in\F$ by $q_i=q^{d_i}$.  
Let $U_q$ be 
the quantized universal enveloping algebra associated to the root datum.
Denote by $U_q^-$ (resp.\ $U_q(\b_-)$) the lower part
(resp.\ the lower Borel subalgebra) of $U_q$.
Let $\{y_1,\ldots,y_M\}$ be a $\Z$-free basis of $d^{-1}Y$.

Define the algebra $\B_q$ 
to be the associative algebra over $\F$
generated by $f_i$ ($i\in I$) and $q^\lambda$ ($\lambda\in d^{-1}Y$)
with defining relations:
\begin{align*}
 &
 \text{$U_q^0$ is a subalgebra of $\B_q$};
 \\ &
 q^\lambda f_i q^{-\lambda} = q^{-\bra\lambda,\alpha_i\ket} f_i 
 \quad\text{for $i\in I$, $\lambda\in d^{-1}Y$};
 \\ &
 f_if_j = q_i^{-\epsilon_{ij}a_{ij}} f_jf_i
 \quad\text{for $i,j\in I$}.
\end{align*}
The last relations are sufficient conditions of 
$q$-Serre relations for $\{f_i\}_{i\in I}$
and called {\em truncated $q$-Serre relations}.
Sending $F_i$ to $f_i$ for each $i\in I$,
we can regard $\B_q$ as a quotient algebra of $U_q(\b_-)$.
Denote the image of $K_i=q_i^{\av_i}$ in $\B_q$ by $k_i$.

Then $\B_q$ is a subalgebra of a $q$-Laurent polynomial ring
over $\F$ generated by $f_i$ ($i\in I$) and $q^{y_\mu}$ ($\mu=1,\ldots,M$)
and hence a tempered domain.
By\lemmaref{lemma:tensor-tempered}, 
$\B_q^{\otimes N}$ is also a tempered domain for any positive integer $N$.

We can define the algebra homomorphism $\Delta_N:U_q\to U_q^{\otimes N}$ 
by $\Delta_1=\id_{U_q}=1$, $\Delta_\nu=(\Delta_{\nu-1}\o1)\circ\Delta$
for $\nu=2,\ldots,N$. 
Then we have $\Delta_N(F_i)=\sum_{\nu=1} F_{i\nu}$,
where $F_{i\nu}=(K_i^{-1})^{\o(\nu-1)}\otimes F_i\o1^{\o(N-\nu)}$.
Therefore we can define the algebra homomorphism 
$\phi_{q,N}:U_q^-\to\B_q^{\otimes N}$
by $\phi_{q,N}(F_i)=\sum_{\nu=1}^N f_{i\nu}$,
where $f_{i\nu}=(k_i^{-1})^{\o(\nu-1)}\otimes f_i\o1^{\o(N-\nu)}$.
Denote the image of $\phi_{q,N}$ by $\N_{q,N}$.
Then $\N_{q,N}$ is also a tempered domain and hence an Ore domain.
Denote $\N_{q,1}$ by $\N_q$ for short.


\section{Non-integral powers}
\label{sec:justifying-non-integral-powers}


\subsection{Evaluation mapping between fields of fractions}
\label{sec:eval}

Let $\F$ be any base field.
Let $A$ be a tempered domain over $\F$.
Set $\A = A\otimes\F[x_1,\ldots,x_M] = A[x_1,\ldots,x_M]$.
Then $\A$ is also a tempered domain over $\F$.
Denote by $K$ the field of fractions of $A$ and by $\K$ that of $\A$.
Using the universality of $K$, we can regard $K$ as a subfield of $\K$.
$\A' = A[x_1^{\pm1},\ldots,x_M^{\pm1}]$ is also a tempered domain.
We can identifies the field of fractions of $\A'$ with $\K$.

For $c=(c_1,\ldots,c_M)\in\F^M$, 
we define the evaluation algebra homomorphism $\ev_c:\A\to A$ at $c$ 
by $\ev_c(f)=f(c)=f(c_1,\ldots,c_M)\in A$ for $f\in\A=A[x_1,\ldots,x_M]$.

Any element $f$ of $\K$ can be represented as $f=g^{-1}h$
for some $g,h\in\A$ with $g\ne 0$.
Fix $c\in\F^M$.
Assume that $g,h,g',h'\in\A$, $g(c),g'(c)\ne0$, 
and $g^{-1}h=g'^{-1}h'$.
Since $\A$ is an Ore domain, 
there exist non-zero $u,u'\in\A$ such that $ug=u'g'$.
Then $uh=ugg^{-1}h=ug'g'^{-1}h'=u'h'$.
We have $u(c)g(c)=u'(c)g'(u)$ and $u(c)h(c)=u'(c)h'(c)$.
Therefore \(
g(c)^{-1}h(c)=(u(c)g(c))^{-1}u(c)h(c)
=(u'(c)g'(c))^{-1}u'(c)h'(c)=g'(c)^{-1}h'(c)
\) in $K$.
This means that if $f\in\K$ can be represented as $f=g^{-1}h$ 
for some $g,h\in\A$ with $g(c)\ne 0$, 
then $\ev_c(f)=f(c)=g(c)^{-1}h(c)\in K$ 
is well-defined and 
does not depend on the choice of $g$ and $h$.

Let $C$ be a subset of $\F^M$ with 
the following {\em Zariski dense property}:
\begin{itemize}
 \item[(D)]
  For every $a\in \F[x_1,\ldots,x_M]$, 
  if $a(c)=0$ for all $c\in C$, then $a=0$ in $\F[x_1,\ldots,x_M]$.
\end{itemize}
For example, for any infinite subset $C_1$ of $\F$, 
the direct product $C_1^M\subset\F^M$ has the property (D).
For every $f\in\A=A\otimes\F[x_1,\ldots,x_M]$, 
if $\ev_c(f)=f(c)=0$ for all $c\in C$, then $f=0$ in $\A$.
Immediately we obtain the following result.

\begin{lemma}
\label{lemma:(D)-1}
 Let $g,h\in\A$, $g\ne 0$, and $f=g^{-1}h$.
 Assume that there exit a subset $C$ of $\F^M$ with the property (D)
 such that $g(c)\ne 0$ for all $c\in C$.
 If $\ev_c(f)=f(c)=0$ for all $c\in C$, then $f=0$ in $\K$.
 \BOX
\end{lemma}

Let $C$ be a subset of $\F^M$ with the property (D).
Then, for any non-zero $b\in \F[x_1,\ldots,x_M]$,
the subset $C_{b\ne 0}=\{\,c\in C\mid b(c)\ne 0\,\}$ of $C$ 
also has the property (D).
In fact, for every $a\in \F[x_1,\ldots,x_M]$, 
if $a(c)=0$ for all $c\in C_{b\ne 0}$, 
then $a(c)b(c)=0$ for all $c\in C$.
Therefore $ab=0$ and hence $a=0$ in $\F[x_1,\ldots,x_M]$.
It follows that, for any non-zero $g_1,\ldots,g_N\in\A$,
the subset \(
 C_{g_1,\ldots,g_N\ne 0}
 =\{\,c\in C\mid g_i(c)\ne 0\ \text{for}\ \text{all}\ i=1,\ldots,N\,\}
\) of $C$ also has the property (D).
From this we can obtain the following result.

\begin{lemma}
\label{lemma:(D)-2}
 Let $C$ be a subset of $\F^M$ with the property (D).
 Take any $f,f'\in\K$.
 By the definition of the field of fractions $\K$,
 there exist $g,h,g',h',g'',h'',u,u'\in\A$
 such that $g,g',g'',u\ne 0$, $f=g^{-1}h$, $f'=g'^{-1}h'$,
 $ff'=(g''g)^{-1}h''h'$, and $f+f'=(ug)^{-1}(uh+u'h')$.
 Then the subset $C'=C_{g,g',g'',u\ne0}$ of $C$ satisfies the property (D)
 and that
 $\ev_c(f)\ev_c(f')=\ev_c(ff')$ and $\ev_c(f)+\ev_c(f')=\ev_c(f+f')$
 for all $c\in C'$. 
 \BOX
\end{lemma}

We shall use these results to justify the conjugation actions of 
non-integral powers in \secref{sec:non-integral-powers}.


\subsection{Non-integral powers in fields of fractions}
\label{sec:non-integral-powers}

In this subsection, we shall justify the conjugation actions of 
non-integral powers 
along the lines of the work \cite{iohara-malikov-GKconj}
by Iohara and Malikov.

Let $\F$ be any base field.
Let $A$ be a tempered domain over $\F$
with generators $f_i$ ($i\in I$)
and defining relations 
$R_\lambda(\{f_i\}_{i\in I})=0$ ($\lambda\in\Lambda$),
where $R_\lambda$ ($\lambda\in\Lambda$) are elements of
the tensor algebra $T(V)$ of $V=\bigoplus_{i\in I}\F f_i$.
That is, $A$ is the quotient algebra of $T(V)$ modulo 
the two-sided ideal generated by $\{R_\lambda\}_{\lambda\in\Lambda}$.
The polynomial ring $A[x]=A\otimes\F[x]$ of one variable over $\A$
is also a tempered domain over $\F$.
We can regard the field of fractions $Q(A)$ as a subfield of $Q(A[x])$.

Assume that a non-zero element $g$ in $A$,
a countable family $\{c_n\}_{n=0}^\infty$ 
of mutually distinct elements in $\F$, and 
an infinite subset $\Gamma$ of $\Z_{\geqq0}$
satisfy the following condition:
\begin{itemize}
 \item[($*$)]
  For any $i\in I$, there exists 
  $\phi_i\in Q(A[x])$ such that, for all $n\in\Gamma$,
  $\ev_{c_n}(\phi_i)=\phi_i(c_n)\in Q(A)$ is well-defined and
  $g^n f_i g^{-n} = \phi_i(c_n)$.
\end{itemize}
For any $\lambda\in\Lambda$ and $k\in\Gamma$,
we have \(
   R_\lambda(\{\phi_i(c_n)\}_{i\in I})
 \ab=\ab R_\lambda(\{g^n f_i g^{-n}\}_{i\in I})
 \ab=\ab g^k R_\lambda(\{f_i\}_{i\in I})g^{-k}
 \ab=\ab 0
\) in $Q(A)$.
By \lemmaref{lemma:(D)-1}, 
$R_\lambda(\{\phi_i(x_1)\}_{i\in I})=0$ in $Q(A[x])$.
Therefore we can define the algebra homomorphism 
$\conj_{g,x}:A[x]\to Q(A[x])$ 
by $\conj_{g,x}(f_i) = \phi_i(x)$ for $i\in I$ and $\conj_{g,x}(x)=x$.
Using the universality of the field of fractions $Q(A[x])$,
we can extend $\conj_{g,x}$ to the algebra automorphism of $Q(A[x])$.
We call $\conj_{g,x}$ the {\em conjugation action of non-integral power}
of $g$ on $Q(A[x])$.

Assume that $\F[x]$ is identified with a subalgebra of 
the polynomial algebra $\F[x_1^{\pm1},\ab \ldots,\ab x_M^{\pm1}]$
of $M$-variables over $\F$.
Then we can also define the algebra automorphism 
$\conj_{g,x}$ of $Q(A[x_1,\ab\ldots,\ab x_M])$
by $\conj_{g,x}(f_i) = \phi_i(x)$ for $i\in I$ 
and $\conj_{g,x}(x_\mu)=x_\mu$ for $\mu=1,\ldots,M$.

If $c_n=n$, then $\conj_{g,x}$ is denoted by $\conj(g^x)$.
If $c_n=q^n$ and $x$ is identified with $q^\lambda$,
then $\conj_{g,x}$ is denoted by $\conj(g^\lambda)$.


\section{Quantized $q$-analogues of birational Weyl group actions}
\label{sec:NY}

In this section, we shall construct quantized $q$-analogues of
the birational Weyl group actions arising from nilpotent Poisson algebras
proposed by Noumi and Yamada \cite{noumi-yamada-birataction2}.

Let $A=[a_{ij}]$ be a symmetrizable GCM symmetrized by $\{d_i\}_{i\in I}$.
Let $\A_{q,0}$ be a quotient tempered domain of 
the lower part $U_q^-$ of a universal enveloping algebra $U_q$ of type $A$.
Then $\A_{q,0}$ is generated by the images $\{f_i\}_{i\in I}$
of the lower Chevalley generators $\{F_i\}_{i\in I}$.
If $A$ is of finite or affine type, 
then any quotient integral domain of $U_a^-$ is a tempered domain.
See \corref{cor:subquot-affine-tempered}.
In order to construct the examples for arbitrary cases, 
see \secref{sec:truncated-q-Serre}.


\subsection{Non-integral power of $f_i$}
\label{sec:f_i^a}

Fix $i\in I$ and assume that $f_i\ne 0$.
\lemmaref{lemma:F_i^nF_j-all} leads to the following formulae:
\begin{equation*}
 f_i^n f_j f_i^{-n}
 = \sum_{k=0}^{-a_{ij}} 
 q_i^{(k+a_{ij})(n-k)}\qbinom{n}{k}_{q_i}\ad(f_i)^k(f_j)f_i^{-k}
 \quad\text{for $n\in\Z$ if $i\ne j$},
\end{equation*}
where $\ad(f_i)^k(f_j)$ denotes 
the image of $\ad(F_i)^k(F_j)\in U_q^-$ in $\A_{q,0}$:
\begin{equation*}
 \ad(f_i)^k(f_j)
 = \sum_{\nu=0}^k
 (-1)^\nu q_i^{\nu(k-1+a_{ij})} \qbinom{k}{\nu}_{q_i}
 f_i^{k-\nu} f_j f_i^\nu.
\end{equation*}
For $j\in I$, define $\phi_{ij}(x)\in Q(A[x])$ by
\begin{equation*}
 \phi_{ij}(x)
 = 
 \begin{cases}
  \displaystyle
  \sum_{k=0}^{-a_{ij}} 
  q_i^{-(k+a_{ij})k}x^{k+a_{ij}}a_{ij;k}(x)\ad(f_i)^k(f_j)f_i^{-k}
  & \text{if $i\ne j$},
  \\
  \displaystyle
  f_i
  & \text{if $i=j$},
 \end{cases}
\end{equation*}
where $a_{ij;k}(x)\in\F[x,x^{-1}]$ for $k\in\Z_{\geqq0}$ are given by
\begin{equation*}
 a_{ij;k}(x)
 = \frac{[x;0]_{q_i}[x;-1]_{q_i}\cdots[x;-k+1]_{q_i}}{[k]_{q_i}!},
 \quad
 [x;\nu]_{q_i} = \frac{xq_i^\nu-x^{-1}q_i^{-\nu}}{q_i-q_i^{-1}}.
\end{equation*}
Then there exist $g_{ij}\in\A_{q,0}$ and $h_{ij}(x)\in\A_{q,0}[x,x^{-1}]$ 
such that $\phi_{ij}(x)=g_{ij}^{-1}h_{ij}(x)$.
Therefore $\phi_{ij}(q^n)$ is well-defined 
and $f_i^n f_j f_i^{-n} = \phi_{ij}(q_i^n)$ 
for all $j\in I$ and $n\in\Z$.

Let 
$(Y,\ab X,\ab\bra\,,\ket,\ab\{\av_i\}_{i\in I},\ab\{\alpha_i\}_{i\in I})$ 
be a root datum of type $A$
and $\{y_1,\ldots,y_M\}$ a $\Z$-free basis of $d^{-1}Y$.
Let $\A_q$ be the tensor product algebra 
$\A_{q,0}\otimes U_q^0=\A_{q,0}[q^{\pm y_1},\ldots,q^{\pm y_M}]$.
Note that $q^\lambda$ ($\lambda\in d^{-1}Y$) 
commute $f_j$ ($j\in I$) in $\A_q$.
Take any $\lambda\in Y$.
Identifying $x$ with $q_i^{\lambda}=q^{d_i\lambda}$,
we regard $\F[x]$ as a subalgebra of $U_q^0$.
Using the result of \secref{sec:non-integral-powers}, 
we can define the algebra automorphism 
$\conj(f_i^\lambda)$ of $Q(\A_q)$
by $\conj(f_i^\lambda)(f_j) = \phi_{ij}(q_i^\lambda)$ for $j\in I$ 
and $\conj(f_i^\lambda)(q^\mu)=q^\mu$ for $\mu\in d^{-1}Y$.
More explicitly we have
\begin{equation}
 \label{eq:q-f^a}
 \conj(f_i^\lambda)(f_j) =
 \begin{cases}
  \displaystyle
  \sum_{k=0}^{-a_{ij}} 
  q_i^{(k+a_{ij})(\lambda-k)}
  \qbinom{\lambda}{k}_{q_i} \ad(f_i)^k(f_j)f_i^{-k}
  & \text{if $i\ne j$},
  \\[\medskipamount]
  \displaystyle
  f_i
  & \text{if $i=j$}.
 \end{cases}
\end{equation}
Note that the right-hand side 
is a Laurent polynomial in $q_i^\lambda$.

For the $q=1$ cases, we have the same construction as the above.
Let $\A_0$ be a quotient tempered domain of 
the universal enveloping algebra of
the lower part $\n_-$ of a Kac-Moody algebra $\g$ of type $A$.
Then $\A_0$ is generated by the images $\{f_i\}_{i\in I}$
of the lower Chevalley generators $\{F_i\}_{i\in I}$.
Let $\A$ be the tensor product algebra 
$\A_0\otimes U(\h)=\A_0[y_1,\ldots,y_\mu]$.
Then we can define the algebra automorphism $\conj(f_i^\lambda)$ 
of $Q(\A)$ by $\conj(f_i^\lambda)(h)=h$ for $h\in\h$ and
the $q\to 1$ limit of \eqref{eq:q-f^a}:
\begin{equation}
 \label{eq:f^a}
 \conj(f_i^\lambda)(f_j)=
 \begin{cases}
  \displaystyle
  \sum_{k=0}^{-a_{ij}}\binom{\lambda}{k}\ad(f_i)^k(f_j)f_i^{-k}
  & \text{if $i\ne j$},
  \\
  \displaystyle
  f_i
  & \text{if $i=j$},
 \end{cases}
\end{equation}
where $\ad(X)(Y)=[X,Y]$ 
and $\binom{\lambda}{k}=\lambda(\lambda-1)\cdots(\lambda-k+1)/k!$.
The left-hand side of \eqref{eq:f^a} is a polynomial in $\lambda$.

\begin{remark}
\label{remark:quantization-of-NY}
\normalfont
 Formula \eqref{eq:f^a} (resp.\ \eqref{eq:q-f^a})
 can be regarded as a quantization (resp.\ quantized $q$-analogue) 
 of Formula (1.9) in \cite{noumi-yamada-birataction2}
 proposed by Noumi and Yamada.

 Simply Formula \eqref{eq:f^a} is 
 the $q\to1$ limit of Formula \eqref{eq:q-f^a}.
 Note that the $q\to1$ limit is not a classical limit
 because \eqref{eq:f^a} is a formula in a non-commutative algebra.

 Let us explain how to obtain 
 Formula (1.9) in \cite{noumi-yamada-birataction2}
 as the classical limit of Formula \eqref{eq:f^a}.
 We replace $f_i$ by $\hbar^{-1}\varphi_i$ and $\lambda$ by $\hbar^{-1}\lambda_i$
 and define $\ad_\hbar$ by $\ad_\hbar(X)(Y)=\hbar^{-1}[X,Y]$,
 where $\hbar$ denotes the Planck constant.
 Assume that $i\ne j$.
 Then Formula \eqref{eq:f^a} is equivalent to
 \begin{equation}
  \conj(\varphi_i^{\hbar^{-1}\lambda_i})(\varphi_j)
  = \sum_{k=0}^{-a_{ij}}
    \frac{\lambda_i(\lambda_i-\hbar)\cdots(\lambda_i-(k-1)\hbar)}{k!}
    \ad_\hbar(\varphi_i)^k(\varphi_j)\varphi_i^{-k}.
  \label{eq:hbar-f^a}
 \end{equation}
 The classical limit of $\hbar^{-1}[X,Y]$ 
 should be the Poisson bracket $\{X,Y\}$.
 Thus, as the classical limit of \eqref{eq:hbar-f^a}, 
 we can obtain Formula (1.9) in \cite{noumi-yamada-birataction2}:
 \begin{equation*}
  s_i(\varphi_j) 
  = \sum_{k=0}^{-a_{ij}} 
    \frac{1}{k!}\left(\frac{\lambda_i}{\varphi_i}\right)^k
    \ad_{\{\,\}}(\varphi_i)^k(\varphi_j),
  \quad
  \ad_{\{\,\}}(X)(Y) = \{X,Y\}.
 \BOX
 \end{equation*}
\end{remark}


\subsection{Quantization of birational Weyl group actions}
\label{sec:quantization-NY}

In the previous subsection, we have constructed the conjugation action 
$\conj(f_i^\lambda)$ of a non-integral powers $f_i^\lambda$ 
on the field of fractions $Q(\A_q)$,
where $i\in I$, $\lambda\in Y$, and 
$\A_q$ is the tensor product algebra of 
a quotient tempered domain $\A_{q,0}$ of $U_q^-$ 
and the Cartan subalgebra $U_q^0$ of $U_q$.
We denote by $f_i$ the image of $F_i$ in $\A_{q,0}$.
We identify $\A_q$ with the Laurent polynomial ring 
$\A_{q,0}[q^{\pm y_1},\ldots,q^{\pm y_M}]$, 
where $\{y_1,\ldots,y_M\}$ is a $\Z$-free basis of $d^{-1}Y$.
Note that $q^\lambda$ ($\lambda\in d^{-1}Y$) commute $f_i$ in $\A_q$.

The Weyl group $W=\bra s_i|i\in I\ket$ 
acts on $Y$. (See \secref{sec:braid-Weyl}.)
This action naturally extends to those on 
$d^{-1}Y$ and $U_q^0=\bigoplus_{\lambda\in d^{-1}Y}\F q^\lambda$.
In this subsection, we denote by $\tw$ 
the action of $w\in W$ on $U_q^0$ 
regarded as a subalgebra of $\A_q$:
\begin{equation*}
 \ts_i(q^\lambda) = q^{s_i(\lambda)}, \quad
 s_i(\lambda) = \lambda - \bra\lambda,\alpha_i\ket\av_i
 \quad \text{for $i\in I$, $\lambda\in d^{-1}Y$}.
\end{equation*}
The action $\tw$ of $w\in W$ on $U_q^0$ is extended to the action on $\A_q$
by $\tw(f_i)=f_i$ for $i\in I$.
The induced action of $\tw$ on $Q(\A_q)$ is also denoted by $\tw$.

\begin{lemma}
\label{lemma:f^a-s_i}
 For any $i,j\in I$ and $\lambda\in Y$, \(
    \conj(f_j^{\lambda}) \circ \ts_i
  = \ts_i \circ \conj(f_j^{s_i(\lambda)})
 \) on $Q(\A_q)$.
\end{lemma}

\begin{proof}
 Take any $k\in I$ and $\mu\in d^{-1}Y$. 
 Then we have
 \begin{align*}
  &
  \conj(f_j^{\lambda}) \circ \ts_i(f_k)
  = \conj(f_j^{\lambda})(f_k)
  = \varphi_{jk}(q_j^\lambda),
  \\ &
  \ts_i \circ \conj(f_j^{s_i(\lambda)})(f_k)
  = \ts_i(\varphi_{jk}(q_j^{s_i(\lambda)}))
  = \varphi_{jk}(q_j^{s_i^2(\lambda)})
  = \varphi_{jk}(q_j^\lambda),
  \\ &
  \conj(f_j^{\lambda}) \circ \ts_i(q^\mu)
  = \conj(f_j^{\lambda})(q^{s_i(\mu)})
  = q^{s_i(\mu)},
  \\ &
  \ts_i \circ \conj(f_j^{s_i(\lambda)})(q^\mu)
  = \ts_i(q^\mu)
  = q^{s_i(\mu)}.
 \end{align*}
 This proves the above lemma.
\end{proof}

\begin{theorem}[quantized birational Weyl group action]
\label{theorem:quantized-NY}
 Assume that $f_i\ne 0$ for all $i\in I$. 
 For each $i\in I$ we define the algebra automorphism $S_i$ of $Q(\A_q)$ by
 \( 
  S_i 
  = \ts_i\circ\conj(f_i^{-\av_i})
  = \conj(f_i^{\av_i})\circ\ts_i.
 \) 
 Then the action of the Weyl group $W$ on $Q(\A_q)$ 
 is defined by $s_i(x)=S_i(x)$ for $i\in I$ and $x\in Q(\A_q)$.
 Explicitly, the following formulae define a representation
 of the Weyl group in algebra automorphisms of $Q(\A_q)$:
 \begin{align*}
  &
  s_i(f_j) = 
  \begin{cases}
   \displaystyle
   \sum_{k=0}^{-a_{ij}} 
   q_i^{(k+a_{ij})(\av_i-k)}
   \qbinom{\av_i}{k}_{q_i} \ad(f_i)^k(f_j)f_i^{-k}
   & \text{if $i\ne j$},
   \\[\medskipamount]
   \displaystyle
   f_i
   & \text{if $i= j$},
  \end{cases}
  \\ &
  s_i(q^\lambda) = q^{s_i(\lambda)} 
  = q^{\lambda-\bra\lambda,\alpha_i\ket\av_i}
  \quad\text{$\lambda\in d^{-1}Y$}.
 \end{align*}
\end{theorem}

\begin{proof}
 It is sufficient to show 
 the braid relations of $\{S_i\}_{i\in I}$ and $S_i^2=1$ for $i\in I$.

 First let us prove $S_i^2=1$. 
 It is sufficient to show that $S_i^2(f_j)=f_j$.
 Since $S_i^2(f_i)=f_i$ is trivial, we can assume $i\ne j$.
 Using \lemmaref{lemma:f^a-s_i} we have 
 $S_i^2=\conj(f_i^{\av_i})\circ\conj(f_i^{-\av_i})$.
 $S_i^2(f_j)$ is a Laurent polynomial $\Phi(q_i^{\av_i})$ 
 of $q_i^{\av_i}$ with coefficients in $Q(\A_q)$.
 Then we have $\Phi(q_i^n)=f_i^{n}(f_i^{-n} f_j f^n)f^{-n}=f_j$ 
 for all $n\in\Z$. 
 Therefore we obtain $S_i^2(f_j)=\Phi(q_i^{\av_i})=f_j$.

 Second let us prove the braid relations for $\{S_i\}_{i\in I}$.
 Assume that $i\ne j$ 
 and $(a_{ij},a_{ji})=(0,0),(-1,-1),(-1,-2),$ or $(-1,-3)$.
 We define the sequences $(i_1,\ldots,i_n)$, $(j_1,\ldots,j_n)$ as follows.
 If $(a_{ij},a_{ji})=(0,0)$, then $n=2$, 
 $(i_1,i_2)=(i,j)$, and $(j_1,j_2)=(j,i)$.
 If $(a_{ij},a_{ji})=(-1,-1)$, then $n=3$, 
 $(i_1,i_2,i_3)=(i,j,i)$, and $(j_1,j_2,j_3)=(j,i,j)$.
 If $(a_{ij},a_{ji})=(-1,-2)$, then $n=4$, 
 $(i_1,\ldots,i_4)=(i,j,i,j)$, and $(j_1,\ldots,j_4)=(j,i,j,i)$.
 If $(a_{ij},a_{ji})=(-1,-3)$, then $n=6$, 
 $(i_1,\ldots,i_6)=(i,j,i,j,i,j)$, and $(j_1,\ldots,j_6)=(j,i,j,i,j,i)$.
 Then the braid relation to be shown is written as
 $S_{i_1}\ldots S_{i_n}=S_{j_1}\ldots S_{j_n}$.

 For $p=1,\ldots,n$, we set
 $\lambda_p=s_{i_n}s_{i_{n-1}}\cdots s_{i_{p+1}}(\av_{i_p})$ and
 $\mu_p=s_{j_n}s_{j_{n-1}}\cdots s_{j_{p+1}}(\av_{j_p})$.
 Then $\lambda_p = a_p\av_i+b_p\av_j$ 
 and  $\mu_p     = c_p\av_i+d_p\av_j$ 
 for some $a_p,b_p.c_p,d_p\in\Z_{\geqq0}$.
 (See \exampleref{example:braid-relations}.)
 For $k,l\in\Z_{\geqq 0}$, we set 
 $u_{k,l} = f_{i_1}^{a_1k+b_1l}\cdots f_{i_n}^{a_nk+b_nl}$ and
 $v_{k,l} = f_{j_1}^{c_1k+d_1l}\cdots f_{j_n}^{c_nk+d_nl}$.
 Then $u(k,l)$ (resp.\ $v(k,l)$) is the image in $\A_{q,0}$ of 
 the left-hand side (resp.\ right-hand side) of the corresponding formula
 in \exampleref{example:Verma-relations}.
 For example, 
 if $(a_{ij},a_{ji})=(-1,-1)$,
 then $u(k,l)=f_i^{l}\ab f_j^{k+l}\ab f_i^{k}$
 The Verma relations mean that $u(k,l)=v(k,l)$ for all $k,l\in\Z_{\geqq0}$.

 By \remarkref{remark:Verma-relations} and \lemmaref{lemma:f^a-s_i}, 
 the condition $S_{i_1}\ldots S_{i_n}=S_{j_1}\ldots S_{j_n}$ is
 equivalent to 
 \( 
  \conj(f_{i_1}^{-\lambda_1})
  \cdots\conj(f_{i_n}^{-\lambda_n})
  =
  \conj(f_{j_1}^{-\mu_1})
  \cdots\conj(f_{j_n}^{-\mu_n}).
 \) 
 Denote the left-hand side by $\phi$ and right-hand side by $\psi$.
 Fix any $t\in I$.
 Then $\phi(f_t)$ and $\psi(f_t)$ belong to 
 $Q(\A_{q,0}[q_i^{-\av_i},q_j^{-\av_j}])$.
 We denote $\phi(f_k)$ by $\Phi(q_i^{-\av_i},q_j^{-\av_j})$ 
 and $\psi(f_k)$ by $\Psi(q_i^{-\av_i},q_j^{-\av_j})$.

 From \lemmaref{lemma:(D)-2} and the definition of $\conj(f_i^\lambda)$,
 it follows that there exists a subset $\Gamma$ of $\Z_{\geqq0}^2$ 
 with the following properties:
 \begin{enumerate}
  \item 
    If $f(x,y)\in \A_{q,0}[x,y]$ 
    and $f(q_i^k,q_j^l)=0$ for all $(k,l)\in\Gamma$,
    then $f(x,y)=0$ in $\A_{q,0}[x,y]$.
  \item 
    $\Phi(q_i^k,q_j^l)$ and $\Psi(q_i^k,q_j^l)$ are well-defined
    for all $(k,l)\in\Gamma$.
  \item 
    $\Phi(q_i^k,q_j^l)=u(k,l)f_tu(k,l)^{-1}$ and
    $\Phi(q_i^k,q_j^l)=v(k,l)f_tv(k,l)^{-1}$
    for all $(k,l)\in\Gamma$.
 \end{enumerate}
 Using the Verma relations, we obtain that
 $\Phi(q_i^k,q_j^l)=\Psi(q_i^k,q_j^l)$ for all $(k,l)\in\Gamma$.
 From \lemmaref{lemma:(D)-1} it follows that \(
  \phi(f_t)
  =\Phi(q_i^{-\av_i},q_j^{-\av_j})
  =\Psi(q_i^{-\av_i},q_j^{-\av_j})
  =\psi(f_t)
 \).  We have just completed the proof.
\end{proof}

\begin{remark}
\label{remark:quantized-NY}
\normalfont
 In \theoremref{theorem:quantized-NY},
 we construct the representation of the Weyl group $W$ 
 in algebra automorphisms of $Q(\A_q)$.
 This can be regarded as a quantized $q$-analogue of 
 the birational Weyl group action
 arising from a nilpotent Poisson algebra
 proposed by Noumi and Yamada
 in \cite{noumi-yamada-birataction2}.
 See also \remarkref{remark:quantization-of-NY}.
\BOX
\end{remark}


\section{Quantized birational Weyl group actions of Hasegawa}
\label{sec:KNY-Hasegawa}

In this section, we shall reconstruct the quantized birational Weyl group 
actions of Hasegawa \cite{hasegawa}.

Let $A=\{a_{ij}\}_{i,j\in I}$ be a symmetrizable GCM 
symmetrized by $\{d_i\}_{i\in I}$,
$d$ the least common denominator of $\{d_i\}_{i\in I}$,
$(Y,\ab X,\ab\bra\,,\ket,\ab\{\av_i\}_{i\in I},\ab\{\alpha_i\}_{i\in I})$ 
a root datum of type $A$.
Let $\{y_1,\ldots,y_M\}$ be a $\Z$-free basis of $d^{-1}Y$.
Assume that if $i\ne j$ and $a_{ij}\ne 0$, 
then $\epsilon_{ij}=\pm1$ and $\epsilon_{ji}=-\epsilon_{ij}$, 
otherwise $\epsilon_{ij}=0$.
Set the base field $\F$ by $\F=\Q(q^{1/d})$
and $q_i\in\F$ by $q_i=q^{d_i}$.  

Consider the tempered domain $\B_q$ defined in \secref{sec:truncated-q-Serre}.
For $i\in I$, define $f_{i1},f_{i2}\in\B_q\otimes\B_q$
by $f_{i1}=f_i\otimes 1$ and $f_{i2}=k_i^{-1}\otimes f_i$.
Note that $f_{i1}+f_{i2}$ ($i\in I$) are the images of 
the lower Chevalley generators $F_i$ in $\B_q\otimes\B_q$.
Therefore $f_{i1}+f_{i2}$ ($i\in I$) 
satisfy the $q$-Serre relations (\secref{sec:def-QUEA})
and hence the Verma relations (\secref{sec:Verma-relations}).

Let $\widetilde{\A}_{q,0}$ be the subalgebra of $\B_q\otimes\B_q$
generated by $f_{i1},f_{i2}$ ($i\in I$).
Then $\widetilde{\A}_{q,0}$ is identified with 
the algebra over $\F$ generated by $f_{i1},f_{i2}$ ($i\in I$)
with defining relations:
\begin{alignat*}{2}
 &
 f_{i\nu}f_{j\nu} = q_i^{-\epsilon_{ij}a_{ij}}f_{j\nu}f_{i\nu}
 & \quad & \text{for $i,j\in I$, $\nu=1,2$},
 \\ &
 f_{i2}f_{j1} = q_i^{a_{ij}}f_{j1}f_{i2}
 & \quad & \text{for $i,j\in I$}.
\end{alignat*}
Note that $q_i^{a_{ij}}=q_j^{a_{ji}}$ because $d_ia_{ij}=d_ja_{ji}$.

Let $\widetilde{\A}_q$ be 
the tensor product algebra $\widetilde{\A}_{q,0}\otimes U_q^0$.
Then $\widetilde{\A}_q$ can be identified with
the Laurent polynomial ring 
$\widetilde{\A}_{q,0}[q^{\pm y_1},\ldots,q^{\pm y_M}]$
with coefficients in $\widetilde{\A}_{q,0}$.
For $i\in I$ and $\lambda\in d^{-1}Y$,
we identify $f_{i1}\o1,f_{i2}\o1,1\o1\otimes q^\lambda\in \widetilde{\A}_q$ with 
$f_{i1},f_{i2},q^\lambda\in\widetilde{\A}_{q,0}[q^{\pm y_1},\ldots,q^{\pm y_M}]$
respectively.
Note that $q^\lambda$ commutes $f_{i1}$ and $f_{i2}$ in $\widetilde{\A}_q$
for $\lambda\in d^{-1}Y$ and $i\in I$.
Since $\widetilde{\A}_q$ is also a tempered domain and hence an Ore domain,
there exists the field of 
fractions $Q(\widetilde{\A}_q)$ of $\widetilde{\A}_q$.

For $i\in I$, define $g_i\in Q(\widetilde{\A}_q)$
by $g_i=f_{i1}^{-1}f_{i2}$.
Let $\A_{q,0}$ (resp.\ $\A_q$) be the subalgebra of $Q(\widetilde{\A}_q)$
generated by $g_i$ ($i\in I$) 
(resp.\ generated by $g_i$ ($i\in I$) and $q^\lambda$ ($\lambda\in d^{-1}Y$)).
Then $\A_{q,0}$ can be identified with the algebra over $\F$ 
generated by $g_i$ ($i\in I$) with defining relations:
\begin{equation}
 g_i g_j = q_i^{-2\epsilon_{ij}a_{ij}} g_j g_i
 \quad \text{for $i,j\in I$}.
 \label{eq:gg=qgg}
\end{equation}
Furthermore $\A_q$ can be identified with 
the Laurent polynomial ring $\A_{q,0}[q^{\pm y_1},\ldots,q^{\pm y_M}]$.
Note that $q^\lambda$ commutes $g_i$ in $\A_q$ for $\lambda\in d^{-1}Y$ and $i\in I.$
Since $\A_{q,0}$ and $\A_q$ are tempered domains,
there exist the fields of fractions $Q(\A_{q,0})$ and $Q(\A_q)$.
We have also $Q(\A_{q,0})\subset Q(\A_q)\subset Q(\widetilde{\A}_q)$.

In \cite{hasegawa}, Hasegawa constructed a representation of 
the Weyl group $W=W(A)$ in algebra automorphisms of $Q(\A_q)$.
Our aim is to reconstruct it by the same method as 
in \secref{sec:NY}.


\subsection{Non-integral power of $f_{i1}+f_{i2}$}
\label{sec:(f+f)^a}

Applying the $q$-binomial theorem (\lemmaref{lemma:q-binom}) 
to $f_{i2}f_{i1}=q_i^2 f_{i1}f_{i2}$, we obtain 
\begin{equation*}
 (f_{i1}+f_{i2})^n
 = 
 \frac
 {(q_i^{-2n}g_i)_{q_i,\infty}}
 {(         g_i)_{q_i,\infty}}
 f_{i1}^n
 \in Q(\widetilde{\A}_q)
 \quad \text{for $n\in\Z$},
\end{equation*}
where $(x)_{i,\infty}=\prod_{\nu=0}^\infty(1+q_i^{2\nu}x)$.
The infinite products in the right-hand side
cancel each other out except finite factors.
Using \eqref{eq:gg=qgg}, and 
$f_{i1}g_j = q_i^{(\epsilon_{ij}-1)a_{ij}}g_jf_{i1}$,
we obtain
\begin{align*}
 &
 (f_{i1}+f_{i2})^n g_j (f_{i1}+f_{i2})^{-n}
 \\ &
 =
 q_i^{(\epsilon_{ij}-1)a_{ij}n}
 g_j
 \frac
 {(q_i^{-2\epsilon_{ij}a_{ij}}q_i^{-2n}g_i)_{q_i,\infty}}
 {(q_i^{-2\epsilon_{ij}a_{ij}}         g_i)_{q_i,\infty}}
 \frac
 {(         g_i)_{q_i,\infty}}
 {(q_i^{-2n}g_i)_{q_i,\infty}}
 \\ &
 =
 q_i^{(\epsilon_{ij}-1)a_{ij}n}
 \frac
 {(q_i^{-2n}g_i)_{q_i,\infty}}
 {(         g_i)_{q_i,\infty}}
 \frac
 {(q_i^{2\epsilon_{ij}a_{ij}}         g_i)_{q_i,\infty}}
 {(q_i^{2\epsilon_{ij}a_{ij}}q_i^{-2n}g_i)_{q_i,\infty}}
 g_j
\end{align*}
for $n\in\Z$.  More explicitly we have
$(f_{i1}+f_{i2})^n g_j (f_{i1}+f_{i2})^{-n}=\phi_{ij}(q_i^n)$
for $n\in\Z$, where $\phi_{ij}(x)\in Q(\A_{q,0}[x])$ ($i,j\in I$)
are defined by
\begin{equation*}
 \phi_{ij}(x)
 =
 \begin{cases}
  \displaystyle
  g_j
  \left(
  \prod_{\nu=0}^{-a_{ij}-1}
  \frac
  {1+q_i^{2\nu}      g_i}
  {1+q_i^{2\nu}x^{-2}g_i}
  \right)
  & \text{if $\epsilon_{ij}=+1$},
  \\[\medskipamount]
  \displaystyle
  x^{2(-a_{ij})}
  \left(
  \prod_{\nu=0}^{-a_{ij}-1}
  \frac
  {1+q_i^{2\nu}x^{-2}g_i}
  {1+q_i^{2\nu}      g_i}
  \right)
  g_j
  & \text{if $\epsilon_{ij}=-1$},
  \\[\medskipamount]
  \displaystyle
  x^{-2}g_i
  & \text{if $i=j$},
  \\[\smallskipamount]
  \displaystyle
  g_j
  & \text{if $a_{ij}=0$},
 \end{cases}
\end{equation*}
Take any $\lambda\in Y$.
Identifying $x$ with $q_i^{\lambda}=q^{d_i\lambda}$,
we regard $\F[x]$ as a subalgebra of $U_q^0$.
Using the result of \secref{sec:non-integral-powers}, 
we can define the algebra automorphism 
$\conj((f_{i1}+f_{i2})^\lambda)$ of $Q(\A_q)$
by $\conj((f_{i1}+f_{i2})^\lambda)(g_j) = \phi_{ij}(q_i^\lambda)$ for $j\in I$ 
and $\conj((f_{i1}+f_{i2})^\lambda)(q^\mu)=q^\mu$ for $\mu\in d^{-1}Y$.

\begin{remark}
\label{remark:Q(Nq2xUq0)}
\normalfont
 We have shown that the conjugation action of 
 a non-integral power $(f_{i1}+f_{i2})^\lambda$ 
 on $Q(\A_q)$ is well-defined.
 Recall that the subalgebra of $\B_q\otimes\B_q$ 
 generated by $\{f_{i1}+f_{i2}\}_{i\in I}$ 
 is denoted by $\N_{q,2}$ in \secref{sec:truncated-q-Serre}.
 Although the conjugation action of $(f_{i1}+f_{i2})^\lambda$
 on $Q(\N_{q,2}\otimes U_q^0)$ is well-defined
 by \theoremref{theorem:quantized-NY},
 it does not reconstruct Hasegawa's action.
\BOX
\end{remark}


\subsection{Reconstruction of Hasegawa's Weyl group actions}
\label{sec:reconstruction-hasegawa}

The Weyl group $W=\bra s_i|i\in I\ket$ acts on $U_q^0$.
This is naturally extended to the action on $Q(\A_q)$ 
so that each $w\in W$ trivially acts on $\{g_i\}_{i\in I}$.
In this subsection,
we denote by $\tw$ the action of $w\in W$ on $Q(\A_q)$:
$\tw(g_i)=g_i$ for $i\in I$ 
and $\tw(q^\lambda)=q^{w(\lambda)}$ for $\lambda\in d^{-1}Y$.

\begin{theorem}
\label{theorem:reconstruction-hasegawa}
 For $i\in I$ we define the algebra automorphism $S_i$ of $Q(\A_q)$ by
 \( 
  S_i 
  = \ts_i\circ\conj((f_{i1}+f_{i2})^{-\av_i})
  = \conj((f_{i1}+f_{i2})^{\av_i})\circ\ts_i.
 \) 
 Then the action of the Weyl group $W$ on $Q(\A_q)$ 
 is defined by $s_i(x)=S_i(x)$ for $i\in I$ and $x\in Q(\A_q)$.
 Explicitly, the following formulae define a representation
 of the Weyl group in algebra automorphisms of $Q(\A_q)$:
 \begin{align*}
  &
  s_i(g_j) = 
  \begin{cases}
   \displaystyle
   g_j
   \left(
   \prod_{\nu=0}^{-a_{ij}-1}
   \frac
   {1+q_i^{2\nu}             g_i}
   {1+q_i^{2\nu}q_i^{-2\av_i}g_i}
   \right)
   & \text{if $\epsilon_{ij}=+1$},
   \\[\medskipamount]
   \displaystyle
   q_i^{2(-a_{ij})\av_i}
   \left(
   \prod_{\nu=0}^{-a_{ij}-1}
   \frac
   {1+q_i^{2\nu}q_i^{-2\av_i}g_i}
   {1+q_i^{2\nu}             g_i}
   \right)
   g_j
   & \text{if $\epsilon_{ij}=-1$},
   \\[\medskipamount]
   \displaystyle
   q_i^{-2\av_i}g_i
   & \text{if $i=j$},
   \\[\smallskipamount]
   \displaystyle
   g_j
   & \text{if $a_{ij}=0$},
  \end{cases}
  \\ &
  s_i(q^\lambda) = q^{s_i(\lambda)} 
  = q^{\lambda-\bra\lambda,\alpha_i\ket\av_i}
  \quad\text{$\lambda\in d^{-1}Y$}.
 \end{align*}
\end{theorem}

\begin{proof}
 Since $f_{i1}+f_{i2}$ ($i\in I$) satisfy the Verma relations,
 we can prove the theorem by the same argument as 
 in the proof of \theoremref{theorem:quantized-NY}.
 Appropriately replacing $f_i^{\av_i}$ and $f_i$ 
 in the proof of \theoremref{theorem:quantized-NY}
 by $(f_{i1}+f_{i2})^{\av_i}$ and $g_i$ respectively, 
 we obtain the proof of the above theorem.
\end{proof}

\begin{remark}[Hasegawa's Weyl group action]
\label{remark:reconstruction-hasegawa}
\normalfont
 In \theoremref{theorem:reconstruction-hasegawa},
 we construct the representation of the Weyl group $W$ 
 in algebra automorphisms of $Q(\A_q)$.
 This can be regarded as a reconstruction of 
 the Weyl group action constructed by Hasegawa in \cite{hasegawa}.
 Set $a_i=q_i^{\av_i}$ and $F_i = a_i^{-1}g_i$ for $i\in I$.
 Do not confuse these $F_i$ with the lower Chevalley generators.
 Then $\A_q$ can be identified with the algebra generated by
 $F_i$ ($i\in I$) and $q^\lambda$ ($\lambda\in d^{-1}Y$)
 with defining relations:
 \begin{align*}
  &
  F_i F_j = q_i^{-2\epsilon_{ij}a_{ij}} F_j F_i
  \quad \text{for $i,j\in I$},
  \\ &
  q^0 = 1,\ q^\lambda q^\mu = q^{\lambda+\mu},
  \ q^\lambda F_i = F_i q^\lambda
  \quad\text{for $\lambda,\mu\in d^{-1}Y$}.
 \end{align*}
 The explicit formulae of the Weyl group action 
 in \theoremref{theorem:reconstruction-hasegawa} 
 can be rewritten as
 \begin{align*}
  &
  s_i(F_j) = 
  \begin{cases}
   \displaystyle
   F_j
   \left(
   \prod_{\nu=0}^{-a_{ij}-1}
   \frac
   {1+q_i^{2\nu}a_iF_i}
   {a_i+q_i^{2\nu}F_i}
   \right)
   & \text{if $\epsilon_{ij}=+1$},
   \\[\medskipamount]
   \displaystyle
   \left(
   \prod_{\nu=0}^{-a_{ij}-1}
   \frac
   {a_i+q_i^{2\nu}F_i}
   {1+q_i^{2\nu}a_iF_i}
   \right)
   F_j
   & \text{if $\epsilon_{ij}=-1$},
   \\[\medskipamount]
   \displaystyle
   F_j
   & \text{otherwise},
  \end{cases}
  \\ &
  s_i(q^\lambda)=q^{s_i(\lambda)},
  \quad
  \text{in particular}\ s_i(a_j) = a_j a_i^{-a_{ij}}.
 \end{align*}
 These formulae essentially coincide with those of Hasegawa.
 Compare these with Hasegawa's formulae,
 Equation (8) in \cite{hasegawa} for the $A^{(1)}_l$ case
 and the example for the $B_2$ case below Theorem 4 in \cite{hasegawa}.

 Let us explain the classical limit of Hasegawa's action.
 Set $q=e^\eta$.
 In the above setting, 
 replace $q$ by $q^{\hbar}$ and $\lambda\in d^{-1}Y$ by $\hbar^{-1}\lambda$. 
 Then $q^\lambda$ is replaced by itself and 
 \begin{equation*}
  \hbar^{-1}[F_i,F_j] 
  = \hbar^{-1}(F_i F_j - F_j F_i) 
  \equiv -2\eta\epsilon_{ij}d_ia_{ij}F_iF_j \mod \hbar.
 \end{equation*}
 Hence the classical limit $\A_q^\cl$ of the algebra $\A_q$ is 
 the commutative Poisson algebra generated 
 by $F_i$ ($i\in I$) and $q^\lambda$ ($\lambda\in d^{-1}Y$) 
 with Poisson brackets defined by
 \begin{alignat*}{2}
  &
  \{F_i, F_j\} = -2\eta\epsilon_{ij}d_ia_{ij}F_iF_j
  & \quad &\text{for $i,j\in I$},
  \\ &
  \{q^\lambda,q^\mu\} = \{q^\lambda,F_i\} = 0
  & \quad & \text{for $\lambda,\mu\in d^{-1}Y$, $i\in I$}.
 \end{alignat*}
 The classical limit of the above Weyl group action can be simply written as
 \begin{equation*}
  s_i(F_j) = 
  F_j \left(\frac{1+a_iF_i}{a_i+F_i}\right)^{\eps_{ij}a_{ij}},
  \quad
  s_i(q^\lambda)=q^{s_i(\lambda)},
 \end{equation*}
 where $a_i = q_i^{\av_i}=q^{d_ia_{ij}}$.
 This action preserves the Poisson brackets of $Q(\A_q^\cl)$.
 The classical case of type $A^{(1)}_2$
 was found by Kajiwara, Noumi, and Yamada in \cite{kajiwara-noumi-yamada-qPIV}.
 See Equation (6) of \cite{kajiwara-noumi-yamada-qPIV}.
 In \cite{hasegawa},
 Hasegawa quantized its generalization
 to an arbitrary symmetrizable GCM case.
\BOX
\end{remark}

%


\end{document}